%% file: paper.tex
\documentclass[final,dvipsnames]{siamart220329}

\usepackage{amssymb,amsfonts}
\usepackage{mathtools}

\usepackage{graphics}
\usepackage{graphicx}
\usepackage{float}
\usepackage{tabularx}
\usepackage{booktabs, multirow}
\usepackage{url}
\usepackage{verbatim,upquote}
\usepackage{setspace}
\usepackage{epstopdf}

\usepackage{subcaption}

\definecolor{hotpink}{rgb}{0.9,0,0.5}
\definecolor{darkred}{rgb}{0.7,0,0}
\hypersetup{colorlinks,urlcolor=blue,citecolor=darkred,linkcolor=blue}

\graphicspath{{figs/}}

\newenvironment{codefragment}[1][htb]{
	\begin{algorithm2e}[#1]%
	}{\end{algorithm2e}}

\usepackage[algo2e,ruled,linesnumbered,noline,algosection]{algorithm2e}
\DontPrintSemicolon
\makeatletter
\let\oldnl\nl
\newcommand{\nonl}{\renewcommand{\nl}{\let\nl\oldnl}} 
\renewcommand{\SetKwProg}[4]{
	\algocf@newcmdside@koif{#1}{\KwSty{#2}\ifArgumentEmpty{#2}\relax{\ 
		}\ProgSty{##2}\KwSty{#3}\ifArgumentEmpty{##1}\relax{ 
			##1}\\\SetAlgoVlined##3\SetAlgoNoLine{#4}{##4}}%
	\algocf@newcommand{l#1}{\@ifstar{\csname algocf@l#1star\endcsname}{\csname 
			algocf@l#1\endcsname}}%
	\algocf@newcmdside{algocf@l#1}{3}{\KwSty{#2} 
		\ProgSty{##2}\KwSty{#3}\algocf@bgroupcode\ 
		##3\algocf@egroupcode\@endalgocfline\ifArgumentEmpty{##1}\relax\ 
		{##1}\strut\par}%
	\algocf@newcmdside{algocf@l#1star}{3}{\KwSty{#2} 
		\ProgSty{##2}\KwSty{#3}\algocf@bgroupcode\ ##3\algocf@egroupcode}%
}%

\newcommand{\algorithmstyle}[1]{\renewcommand{\algocf@style}{#1}}
\makeatother
\SetKwProg{Fn}{function}{}{}
\newcommand{\funcomment}[1]{\nonl$\triangleright$ \textit{#1}}
\SetKwComment{Comment}{$\triangleright$\ }{}
\SetCommentSty{itshape}
\renewcommand{\BlankLine}{\vskip 11pt}%
\newcommand{\myhrulefill}{\leavevmode\leaders\hrule height 0.7ex depth 
	\dimexpr0.4pt-0.7ex\hfill\kern0pt}

\newcommand{\algcmd}[1]{\ensuremath{\text{\textbf{\em {#1\ }}}}}

\newcommand{\algand}{\algcmd{and}}

\newcommand{\algto}{\algcmd{to}}
\newcommand{\algdownto}{\algcmd{down\ to}}

\def\DS{\displaystyle}

\def\eu{\ensuremath{\mathrm{e}}}

\renewcommand{\t}[1]{\texttt{#1}}

\def\>{\mskip\medmuskip}

\DeclarePairedDelimiter\norm{\lVert}{\rVert}
\newcommand{\normi}[1]{\ensuremath{\norm{#1}_1}}

\newcommand{\R}{\ensuremath{\mathbb{R}}}
\newcommand{\C}{\ensuremath{\mathbb{C}}}

\newcommand{\nbyn}{\ensuremath{n\times n}}

\newcommand{\floor}[1]{\lfloor #1 \rfloor}
\newcommand{\ceil}[1]{\lceil #1 \rceil}
\DeclarePairedDelimiter\abs{\lvert}{\rvert}

\DeclareMathOperator{\fl}{\operatorname{f\kern.2ptl}} 
\DeclareMathOperator{\op}{op}

\def\pade{Pad\'e }

\def\ps{Paterson--Stockmeyer}
\def\wh{\widehat}

\newcommand{\gn}{\ensuremath{\gamma_n}}

\newcommand{\us}{\ensuremath{u_{s}}}
\newcommand{\uf}{\ensuremath{u_{f}}}
\newcommand{\X}{\ensuremath{\mathcal{X}}}

\mathcode`@="8000 
{\catcode`\@=\active\gdef@{\mkern1mu}}

\makeatletter
\let\oldtabcr\@tabcr

\gdef\@tabcr{\@stopline \@ifstar{\penalty%
		\@M \@xtabcr}\@xtabcr\mynewline}

\newenvironment{code}{%
	\def\mynewline{
		\llap{\footnotesize
			\hspace{22pt}}%
	}
	\mathcode`\:="603A  
	\def\colon{\mathchar"303A}
	\par
	\upshape
	\begin{list} 
		{} {\leftmargin = 1cm}
		\item[]
		\begin{tabbing}
			\hspace*{.3in} \= \hspace*{.3in} \=
			\hspace*{.3in} \= \hspace*{.3in} \= \kill
			\mynewline
		}{\end{tabbing}\end{list}}
	
\newenvironment{dedication}
{
	\small      
	\thispagestyle{empty}
	\vspace*{.1cm}
	\itshape             
	\centering          
}
{\par 
}

\title{Mixed-Precision Paterson--Stockmeyer Method for 
	Evaluating Polynomials of Matrices%
	\thanks{Version of December 4, 2024.}}
\author{Xiaobo Liu%
	\thanks{Department of Mathematics,
		University of Manchester,
		Manchester, M13 9PL, UK
		(\email{xiaobo.liu@manchester.ac.uk}).}
}

\makeatother

\mathchardef\Gamma="7100 \mathchardef\Delta="7101
\mathchardef\Theta="7102 \mathchardef\Lambda="7103
\mathchardef\Xi="7104 \mathchardef\Pi="7105 \mathchardef\Sigma="7106
\mathchardef\Upsilon="7107 \mathchardef\Phi="7108
\mathchardef\Psi="7109 \mathchardef\Omega="710A

\headers{Mixed-Precision Paterson--Stockmeyer Method}%
{Xiaobo Liu}

\begin{document}
\maketitle

\begin{abstract}
The Paterson--Stockmeyer method is an evaluation scheme for matrix polynomials with scalar coefficients that arise in many state-of-the-art algorithms based on polynomial or rational approximation, for example, those for computing transcendental matrix functions. We derive a mixed-precision version of the Paterson--Stockmeyer method that 
is particularly useful for evaluating matrix polynomials with scalar coefficients of decaying magnitude.
The new method is mainly of interest in the arbitrary precision arithmetic, 
and it is attractive for high-precision computations.
The key idea is to perform computations on data of small magnitude in low precision, and rounding error analysis is provided for the use of lower-than-the-working precisions.
We focus on the evaluation of the Taylor approximants of the matrix exponential and show the applicability of our method to the existing scaling and squaring algorithms. We also demonstrate through experiments the general applicability of our method to other problems, such as
computing the polynomials from the Pad\'e approximant of the matrix exponential and the Taylor approximant of the matrix cosine. 
Numerical experiments show our mixed-precision Paterson--Stockmeyer algorithms can be more efficient
than its fixed-precision counterpart while delivering the same level of accuracy. 
\end{abstract}

\begin{keywords}
Paterson--Stockmeyer method, mixed-precision algorithm, rounding error analysis, arbitrary precision arithmetic, polynomial of matrices, matrix function
\end{keywords}

\begin{MSCcodes}
15A16, 65F60, 65G50
\end{MSCcodes}

\begin{dedication}
	This work is dedicated to the memory of 
	Prof. Nick Higham, for the guidance and support he gave the author and for the continuing inspiration that his work brings to the community.
\end{dedication}

\section{Introduction}

The \ps\ (PS) method~\cite{past73} is an evaluation scheme for matrix
polynomials with scalar coefficients that is used in many state-of-the-art
algorithms based on polynomial or rational approximants for computing
transcendental functions of matrices, for example, the matrix
exponential~\cite{fahi19},~\cite{sidr15}, the matrix logarithm~\cite{fahi18}, and the
matrix trigonometric and hyperbolic functions and their
inverses \cite{ahl22}, \cite{ahr15}, \cite{aphi16}, \cite{siap17}.  
In the PS scheme, a matrix polynomial $p_m(X) = \sum_{i=0}^{m}b_iX^i$ at
$X\in\C^{\nbyn}$ is written as
\begin{equation}\label{eq:ps}
p_m(X) = \sum_{i=0}^{r} B_i\cdot(X^s)^i ,\quad r = \floor{m/s},
\end{equation}
where $s$ is an integer parameter, $X^0\equiv I$ denotes the identity matrix of order $n$, and 
\begin{equation*}
B_i = \begin{cases}
\DS\sum_{j=0}^{s-1}b_{s i+j}X^j, & i=0,1,\dots,r-1, \\
\DS\sum_{j=0}^{m-sr} b_{s r+j}X^{j}, & i=r.
\end{cases}
\end{equation*}
The first $s$ positive powers of $X$ are computed once the parameter $s$
is chosen; then~\eqref{eq:ps} is evaluated by the matrix version of 
Horner's method, with each coefficient polynomial $B_i$ formed via 
explicit powers
reusing the computed powers of $X$. 

In the evaluation of polynomials of matrices, matrix multiplications 
have the highest asymptotic cost amongst all the matrix operations, so 
it is sensible to measure the efficiency of an evaluation scheme by 
the number of matrix multiplications required.
This quantity is known to be minimised by setting 
$s=\sqrt{m}$ (which is 
not necessarily an integer) for a given matrix polynomial 
$p_m(X)$, and the practical choices of 
$s=\floor{\sqrt{m}}$ or 
$s=\ceil{\sqrt{m}}$ yield exactly the same 
cost~\cite{fasi19},~\cite[pp.~29--30]{harg05}, which is about 
$2\sqrt{m}$.
The above discussion is most relevant to the case when the polynomial $p$ 
is dense (most of the coefficients $b_i$ are nonzero), and the 
economics of the evaluation can be rather different if $p_m$ is sparse, 
which is not the focus of the work.
One downside of the PS method is that it requires 
$(s+2)n^2=O(\sqrt{m}n^2)$ memory locations including the storage of the 
first $s$ powers of $X$, in contrast to the $n^2$ elements of storage 
from Horner's method.  Van Loan~\cite{vanl79} proposed a modification 
of the PS method which reduces the storage requirement to $4n^2$ by 
computing $p_m$ one column at a time, at the price of about $40\%$ extra 
flops.
A more efficient block variant of Van Loan's algorithm is developed 
in~\cite{hst18}, 
where it is shown that the computation of all the three 
mentioned schemes can be accelerated 
by reducing the argument matrix to 
its Schur form if the degree $m$ is sufficiently large---so the savings from performing matrix multiplications between triangular (instead of full) matrices outweigh the extra costs in reducing the matrix to Schur form.

With the increasing availability of precisions beyond the IEEE double and single precision arithmetics~\cite{ieee08} in both hardware and software,
the landscape of floating-point precisions has been broadening and scientific computing has been carrying out in an intrinsically mixed precision world. 
Utilizing mixed precisions within numerical algorithms can reduce data transfers and memory access in computer processors, which can improve speed and reduce energy consumption, while the precisions need to be 
chosen prudently to maintain the targeted accuracy. 
There has been a steep increase of interest in the study of mixed precision algorithms in numerical linear algebra over the years~\cite{hima22}.
In particular, multiprecision algorithms that use one or more arbitrary precisions have already been developed for
the computation of general analytic matrix functions~\cite{hili21},
the matrix exponential~\cite{fahi19}, the matrix logarithm~\cite{fahi18}, 
and the matrix cosine and its Fr\'echet derivative~\cite{ahl22}.

In this work we aim to utilize multiple precisions 
in the computation of $p$ by the PS method so as to achieve 
$\norm{p-\wh{p}}\lesssim cnu\norm{p}$, where $\wh{p}$ denotes the computed polynomial, $u$ is the unit roundoff of the working precision, and $c$ denotes some mild constant, given that its scalar coefficients $b_i$ enjoy a certain fast decaying property.
Our idea is inspired by the fundamental fact that computations performed on data of small magnitude can use low precision. For example, in the computation of $X=C+AB$ where 
$|A||B|\ll |C|$ then the matrix product $AB$ can be computed in lower precision than the subsequent summation without significantly impacting the overall accuracy. 
The new method will be particularly useful for accelerating the computation of matrix functions via a polynomial or rational function
in arbitrary precision, for which the working precision can be much higher than the IEEE double precision; see~\cite[sect.~2.2.2]{fasi19a},~\cite[sect.~1]{liu22} for the various applications and details of 
languages and libraries that support arbitrary precision algorithms for matrix functions.

We begin in Section~\ref{sect:err-analy} by stating the main theorems which are the building blocks for analysing the errors in the evaluation of $p_m(X)$ in~\eqref{eq:ps} and discuss the evaluation scheme following from the error analysis.
In Section~\ref{sect:exp}, we apply the framework derived in the previous section to Taylor approximants of the matrix exponential and show the applicability of our method to existing scaling and squaring algorithms for that function, particularly when the norm of the input matrix is sufficiently small, in which case accuracy of the mixed-precision method is shown to be guaranteed.
Numerical experiments are presented to demonstrate the accuracy and efficiency of the algorithms.
In Section~\ref{sect.general.framework}, we illustrate with examples the general applicability of our framework to the computation of matrix polynomials with scalar coefficients that decay in modulus.
Conclusions are drawn in Section~\ref{sect.concl}.

Throughout this work we denote by $\norm{\cdot}$ any consistent matrix norm,
by $\mathbb{N^+}$ the set of positive integers, and by
$u$ the unit roundoff of the floating-point arithmetic. 
We adopt the MATLAB-style colon notation for indices, i.e., 
$a\colon b$ represents the index set $\{a, a+1,\dots,b\}$.
An inequality expressed as ``$a\ll b$'' can be read as ``$a$ is sufficiently less than $b$''. In some context we use the term ``fixed precision'' to mean that the floating point precision has been fixed; this is not to be confused with fixed precision in contrast with floating point precision.

\section{Rounding error analysis and evaluation scheme}\label{sect:err-analy}

We use the standard model of floating-point 
arithmetic~\cite[sect.~2.2]{high:ASNA2}
\begin{equation}\label{eq:fpm-scalar}
\fl(x \op y) = (x \op y)(1+\delta),\quad |\delta|\le u, 
\end{equation}
where $x$ and $y$ are floating-point numbers and $\op$ denotes
addition, subtraction, multiplication, or division.
For matrix multiplication, we have \cite[sect.~3.5]{high:ASNA2}
\begin{equation}\label{matmult-err}
	\fl(AB) = AB + E,\quad |E|\le \gamma_n|A||B|,
\end{equation}
where $\gn:= nu/(1-nu)$, assuming $nu<1$. 
If more than one precision is involved in a computation, we will use the operator $\fl_j(\cdot)$ to denote an operation executed in precision $u_j$ and 
$\gamma_n^j:= nu_j/(1-nu_j)$, assuming $nu_j<1$. Define $\theta_{i,j}=u_i/u_{j}$, $i,j=0\colon r$, so we have
\begin{equation}\label{eq:gamma-theta}
	\gamma_{n}^i = \gamma_{n\theta_{i,j}}^{j}.
\end{equation}

The evaluation of $p_m(X)$ in~\eqref{eq:ps} is customarily performed via Horner's method, that is, we compute
\begin{equation}\label{eq:ps-Horner}
	p_m(X) = B_0 + X^s\big(B_1 + X^s\big(B_2 +\cdots+ X^s(B_{r-1} + X^s B_r)\big)\big)
\end{equation}
starting from the quantities in the innermost brackets. The standard fixed-precision Paterson--Stockmeyer scheme for computing polynomials of matrices is presented in Algorithm~\ref{alg.fpps}.

\begin{algorithm2e}
	\setstretch{1}
	\caption{Paterson--Stockmeyer scheme for polynomials of matrices.}
	\label{alg.fpps}
	\SetAlgoVlined
	{\nonl
		\begin{minipage}{0.95\linewidth}
			\pretolerance=100\tolerance=200\hyphenpenalty=5
			Given $X\in\mathbb{C}^{\nbyn}$ and a set of polynomial coefficients $\{b_i\}_{i=0}^m$, this algorithm computes the matrix polynomial with scalar coefficients
			$P\equiv p_m(X)$ in the form of~\eqref{eq:ps} using the Paterson--Stockmeyer scheme in floating-point arithmetic. 
	\end{minipage}}
	\BlankLine
	$s \gets \ceil{\sqrt{m}}$ \; 
	$r \gets \floor{m/s}$ \;
	$\X_0 \gets I$\;
	\For{$i \gets 1$ \algto $s$}{
		$\X_{i} \gets \X_{i-1} X$  \;
	}
	$P = \sum_{j=0}^{m-sr}b_{sr+j}\X_{j}$ \; 
	\For{$i \gets r-1$ \algdownto $0$}{
		$P = \sum_{j=0}^{s-1}b_{si+j}\X_{j} + \X_{s}P$\;
	}
	\Return{$P$}\;
	\SetAlgoNoLine
\end{algorithm2e}

In this paper we are most interested in the case where 
the $\abs{b_i}$ decay quickly, so we have, for some positive integer $\nu\in\left[1,r\right]$,
\begin{equation}\label{eq:Bi-Xs-bnd}
	\norm{B_i}\norm{X^s}= \tau_i\norm{B_{i-1}},\quad \tau_i\ll1,\quad i=\nu\colon r.
\end{equation}
This set of conditions can be written as, for $i=\nu\colon r$,
\begin{equation*}
	\norm{b_{si}I + b_{si+1}X + \cdots +b_{s i+s-1}X^{s-1}}\norm{X^s} \ll \norm{b_{si-s}I + b_{si-s+1}X + \cdots +b_{si-1}X^{s-1}},
\end{equation*}  
which will hold if $\norm{X}$ is sufficiently small. The intuition is that the dominant terms in
$B_i$ and $B_{i-1}$ are linear combination of powers of $X$ generally from the same set with each pair of corresponding scalar coefficients being $s$ indices apart from the series $\{b_i\}$, of which the modulus decays rapidly. 

Consider for example the numerator polynomial of the \pade approximant to the matrix exponential of the $3\times 3$ test matrix from~\cite{ward77},
\begin{equation*}
	A = \begin{bmatrix}
		-131 & 19 & 18 \\
		-390 & 56 & 54 \\
		-387 & 57 & 52
	\end{bmatrix},
\end{equation*}	
for which the state-of-the-art double-precision-oriented algorithm~\cite{alhi09a} (MATLAB \texttt{expm}) performs the scaling $X=A/2^6$ before computing the $[13/13]$ \pade approximant. If the Paterson--Stockmeyer method is used in this case, then $s=\ceil{\sqrt{13}}=4$ and we arrive at
\begin{align*}
\normi{B_1}\normi{X^s} &= 
\Big\|\frac{11}{5520}I+\frac{3}{18400}X + \frac{1}{96600}X^2 + \frac{1}{1932000}X^3\Big\|_1 \normi{X^4}
\approx 8.6\times10^{-3} 
\\ 
&\ll
1.1\times10^{1} \approx 
\Big\|I+\frac{1}{2}X + \frac{3}{25}X^2 + \frac{11}{600}X^3\Big\|_1
=\normi{B_0},
\end{align*}
where the dominant terms in $B_0$ and $B_1$ (such that they approximate $B_0$ and $B_1$, respectively, to the correct order of magnitude) are both from $\{I,X\}$. In this example the conditions~\eqref{eq:Bi-Xs-bnd} hold with $\tau_1\approx 7.8\times 10^{-4}$ and one can verify that $\tau_2$, $\tau_3$ are even smaller.
Later we will discuss to what extent the conditions~\eqref{eq:Bi-Xs-bnd} can hold for the polynomial $p_m(X)$ in general.

Define the polynomial
\begin{equation}\label{eq:ps-Horner-part}
	q(X) := B_{\nu-1} +  X^s\big(B_{\nu} + X^s\big(B_{\nu+1}+\cdots+ X^s(B_{r-1} + X^s B_r)\big)\big),
\end{equation}
which is exactly $p_m(X)$ if $\nu=1$. Assuming that~\eqref{eq:Bi-Xs-bnd} is satisfied, 
our idea for computing $q(X)$ is to start with the lowest precision in forming the matrix product in the innermost brackets, and then gradually and adaptively increase the precision (up to the working precision) for the subsequent matrix products outwards, aiming to still deliver the full working precision accuracy for the computation of $q(X)$.
The following theorem provides a rounding error bound on the forward error of the process, where potentially different precisions $u_{i}$, $i=\nu-1\colon r$ 
are involved,
and the
proof is given in Appendix~\ref{append:proof-ps-mp-err}.  

\begin{theorem}\label{thm:ps-mp-err-norm}
If
	$\norm{\wh{B}_i-B_i}\le u_i\norm{B_i}$, $i=\nu-1\colon r$ and
$\norm{\wh{Y}-Y}\le u_{\nu}\norm{Y}$ where $Y\equiv X^s$, then  for the matrix $\wh{q}$, computed in finite precision, for $q(X)$ in~\eqref{eq:ps-Horner-part}, the evaluation scheme
	\begin{code}
		$\widehat{\varphi}_r = \widehat{B}_r$ \\
		for $j=r:-1:\nu$ \\
		\> $\widehat{\varphi}_{j-1}=\fl_{j-1}(\wh{B}_{j-1} + \fl_j(\wh{Y} \widehat{\varphi}_{j}))$ \\
		end \\
		$\wh{q} = \widehat{\varphi}_{\nu-1}$  	
	\end{code}
	satisfies
	\begin{equation}\label{eq:ps-mp-err-norm}
		\norm{\wh{q} - q(X)} \le 
		\gamma_{f_{r}}^{\nu-1}\norm{X^s}^{r-\nu+1}\norm{B_r} + \gamma_{f_{r-1}}^{\nu-1}\norm{X^s}^{r-\nu}\norm{B_{r-1}} + \cdots + 
		\gamma_{f_{\nu-1}}^{\nu-1}\norm{B_{\nu-1}},
	\end{equation}
where
\begin{align*} 
	f_{r} &= \frac{nu_r}{u_{\nu-1}} + \frac{(n+1)}{u_{\nu-1}}(u_{r-1}+u_{r-2}+\dots+u_{\nu}) + 1,   \\
	f_{r-1} &= \frac{(n+2)u_{r-1}}{u_{\nu-1}} + \frac{(n+1)}{u_{\nu-1}}(u_{r-2}+u_{r-3}+\dots+u_{\nu}) + 1,   \\
	&\vdots  \\
	f_{\nu} &= \frac{(n+2)u_{\nu}}{u_{\nu-1}}+1, \\
	f_{\nu-1} &= 2.
\end{align*}
\end{theorem}

The constants $f_i$ in~\eqref{eq:ps-mp-err-norm} can be bounded above by 
\begin{equation*} 
	f_{i} \le  (n+2)\left(\frac{u_i}{u_{\nu-1}} + \frac{u_{i-1}}{u_{\nu-1}} +\dots+\frac{u_{\nu}}{u_{\nu-1}}\right) + 1
	\le (n+2)
	\sum_{j=\nu-1}^{i}\frac{u_j}{u_{\nu-1}},\quad 
	i=\nu\colon r.
\end{equation*}
If~\eqref{eq:Bi-Xs-bnd} is satisfied, then in Theorem~\ref{thm:ps-mp-err-norm} we choose the precisions
\begin{equation}\label{eq:choose-ui}
	u_i = \frac{\norm{B_{\nu-1}}u_{\nu-1}}{\norm{B_i}\norm{X^s}^{i-\nu+1}},\quad  i=\nu\colon r,
\end{equation}
which implies $u_{i-1}= \tau_i u_i$; and we
take $\tau=\max_i \tau_i\ll 1$ so that we have
\begin{equation}\label{eq:ui-ui-1-relation}
u_{i-1}\le \tau u_i,\quad  i=\nu\colon r.
\end{equation}
It follows 
$u_{i-2} \le \tau^2u_i, u_{i-3} \le \tau^3u_i, \dots, u_{\nu-1} \le \tau^{i-\nu+1}u_i$ and that, for $i=\nu\colon r$,
\begin{align*} 
	\frac{u_{\nu-1}f_{i}}{u_i} &\le 
	(n+2)\left(\frac{u_{i}}{u_i} + \frac{u_{i-1}}{u_i} + 
\dots + \frac{u_{\nu-1}}{u_i}\right) \\
	&\le (n+2)\left(1 + \tau + 
	\dots+\tau^{i-\nu+1}\right) \\
	&= n+2 + O(n\tau),
\end{align*}
and therefore we have $f_i \lesssim \left(n+n\tau+2\right)u_i/u_{\nu-1}$ and thus
\begin{equation*}
	f_i\norm{B_i}\norm{X^s}^{i-\nu+1} \lesssim \frac{(n+n\tau+2)u_i}{u_{\nu-1}}\norm{B_i}\norm{X^s}^{i-\nu+1} =(n+n\tau+2)\norm{B_{\nu-1}}.
\end{equation*}
Then, since $i\gamma_k^{\nu-1}\le \gamma_{ik}^{\nu-1}$~\cite[Lem.~3.3]{high:ASNA2}, we have the bound on the matrix $\wh{q}$ computed via the evaluation scheme in Theorem~\ref{thm:ps-mp-err-norm}:
	\begin{align*}
	\norm{\wh{q} - q(X)} &\le 
	\gamma_{f_{r}}^{\nu-1}\norm{X^s}^{r-\nu+1}\norm{B_r} + \gamma_{f_{r-1}}^{\nu-1}\norm{X^s}^{r-\nu}\norm{B_{r-1}} + \cdots + 
	\gamma_{f_{\nu-1}}^{\nu-1}\norm{B_{\nu-1}} \\
	& \lesssim (r-\nu+1)\gamma_{(n+n\tau+2)\norm{B_{\nu-1}}}^{\nu-1}  + \gamma_{2}^{\nu-1}\norm{B_{\nu-1}} \\
	& \lesssim \frac{(r-\nu+1)(n+n\tau+2)+2}{1-(n+n\tau+2)\norm{B_{\nu-1}}u_{\nu-1}} \norm{B_{\nu-1}}u_{\nu-1} \\
	& \approx \frac{(r-\nu+1)n}{1-(n+n\tau+2)\norm{B_{\nu-1}}u_{\nu-1}} \norm{q(X)}u_{\nu-1}.
\end{align*}
Therefore, if $\left((1+\tau)n+2\right)\norm{B_{\nu-1}}u_{\nu-1}\ll 1$, then we can choose the precisions $u_i$ by~\eqref{eq:choose-ui} 
such that the computed matrix $\wh{q}$ of $q(X)$ has approximately a normwise relative error of $(r-\nu+1)nu_{\nu-1}$, where $r=\floor{m/s}$, and,
in particular, if~\eqref{eq:Bi-Xs-bnd} holds for $\nu=1$ 
then
the computed matrix $\wh{q}_m$ of $q_m(X)$ has approximately a normwise relative error of $rnu_{0}$.

For $s=1$ the PS scheme~\eqref{eq:ps-Horner} reduces to Horner's method, in which case 
the conditions~\eqref{eq:Bi-Xs-bnd} become
\begin{equation}\label{eq:bi-X-bnd}
	|b_{i}|\norm{X}=\tau_i |b_{i-1}|, \quad \tau_i\ll1,\quad i=\nu\colon m,
\end{equation}
for some positive integer $\nu\in\left[1,m\right]$.
We can obtain an analogous result to Theorem~\ref{thm:ps-mp-err-norm}. In this case it can be shown that
if $\left((1+\tau)n+2\right)\abs{b_{\nu-1}}u_{\nu-1}\ll 1$, then we can choose the precisions by
\begin{equation*}
	u_i = \frac{\abs{b_{\nu-1}}u_{\nu-1}}{\abs{b_i}\norm{X}^{i-\nu+1}},\quad  i=\nu\colon m
\end{equation*}
such that the computed matrix $\wh{q}$ of $q(X)$ has approximately a normwise relative error of $(m-\nu+1)nu_{\nu-1}$,
and, in particular, if~\eqref{eq:bi-X-bnd} holds for $\nu=1$ then 
the computed matrix $\wh{q}_m$ of $q_m(X)$ has approximately a normwise relative error of $mnu_{0}$.

The requirement~\eqref{eq:bi-X-bnd} is made between any two consecutive coefficients and it can only hold if $\norm{X}$ is sufficiently small and the decay rate of $|b_i|$ is sufficiently large.
On the other hand, the PS scheme with sufficiently large $s$ can mitigate this potentially very stringent requirement. 
The requirement in~\eqref{eq:bi-X-bnd} is on adjacent coefficients $b_{i-1}$ and $b_{i}$, but in~\eqref{eq:Bi-Xs-bnd} the dominant terms are $s$ indices apart so the condition is more likely to be satisfied.
Also, the error bound associated with the PS scheme is smaller than that of Horner's method by at most a factor of approximately $r/m=\floor{m/s}/m\approx 1/s$. 
We henceforth focus on the PS evaluation scheme~\eqref{eq:ps-Horner} in the general case ($s$ is not necessarily equal to $1$). 

To summarize, the framework of the mixed-precision PS scheme for the matrix polynomial in~\eqref{eq:ps} is that we exploit lower precisions $u_r\ge u_{r-1}\ge \cdots \ge u_{\nu}$ in the computation of $q(X)$ and then perform the matrix products and sums in the evaluation of
\begin{equation}\label{eq:pm-prec-u-part}
	p_m(X) = B_0 + X^s\big(B_1 + X^s\big(B_2 + \cdots+ X^s\big(B_{\nu-2} + X^s q(X)\big)\big)\big)
\end{equation}
in the working precision $u:=u_{\nu-1}$. 
The required powers of $X$ are formed explicitly and each $B_i$ is formed by reusing these computed powers, which involves only matrix scaling and additions,
so we will form the powers of $X$ in the working precision $u$. 
From the earlier discussion, the computed matrix $\wh{q}$ of $q(X)$ has approximately a normwise relative error bounded above by $(r-\nu+1)nu$, which is satisfactory for the evaluation of~\eqref{eq:pm-prec-u-part} in precision $u$. 
However, this error bound is from Theorem~\ref{thm:ps-mp-err-norm} and is subject to 
$\norm{\wh{Y}-Y}\le u_\nu\norm{Y}$, where $Y\equiv X^s$, and $\norm{\wh{B}_i-B_i}\le u_i\norm{B_i}$, $i=\nu-1\colon r$. 
Among the latter requirements, we just need to ensure that  $\norm{\wh{B}_{\nu-1}-B_{\nu-1}}\le u\norm{B_{\nu-1}}$ is satisfied by the choice of the precisions~\eqref{eq:choose-ui}, as we can do the matrix scaling and additions
required in assembling the $B_i$ in the working precision $u$.   
In fact, since $B_{\nu-1}$ is only involved in the final matrix summation in the evaluation of $q(X)$ via~\eqref{eq:ps-Horner-part}, it is not hard to see 
that we can ease the condition to
\begin{equation}\label{eq:cond-B0}
\norm{\wh{B}_{\nu-1}-B_{\nu-1}}\lesssim cnu\norm{B_{\nu-1}},
\end{equation}
where $c$ denotes some mild (with respect to $n$) constant comparable to $r-\nu+1$, still achieving the same error bound on the computed $\wh{q}$.

We next provide rounding errors analysis for the computation of matrix powers and polynomials and discuss its practical implications for our use case.

\subsection{Powers and polynomials of matrices}
We first derive an upper bound on the rounding errors in the computation of
$X^k$ in a fixed precision.  In our results, when we write $X\in\R^{\nbyn}$ it is understood
that $X$ is a matrix of floating-point numbers.


\begin{lemma}\label{thm:xpower-err}
	For $X\in\R^{\nbyn}$ and $X_t = X^t$ with $t\in\mathbb{N}^+$,
	the computed matrix $\wh{X}_t = \fl(X^t)$ 
	obtained by repeated multiplication
	$\wh{X}_k = \fl(\wh{X}_{k-1}X)$, $k=1\colon t$, $\wh{X}_0=X$ in precision $u$ satisfies
	$|\wh{X}_{t} - X^{t}| \le \gamma_{(t-1)n}|X|^{t}$.
	
	\proof 
	For $t = 2$ the bound holds by \eqref{matmult-err}.
	Suppose the bound holds for $t=k-1$:
	\begin{equation*}
		|\wh{X}_{k-1} - X^{k-1}| \le \gamma_{(k-2)n}|X|^{k-1}.
	\end{equation*}
	We have, for $t=k$,
	\begin{align*}
		|\wh{X}_{k} - X^{k}| &\le |\wh{X}_{k} - \wh{X}_{k-1}X| 
		+ |\wh{X}_{k-1}X - X^{k}| 
		\le\gamma_{n}|\wh{X}_{k-1}||X| + |\wh{X}_{k-1} - X^{k-1}||X|\\
		&\le \gamma_{n}(1+\gamma_{(k-2)n})|X|^{k-1}|X| + 
		\gamma_{(k-2)n}|X|^{k} 
		\le \gamma_{(k-1)n}|X|^{k},
	\end{align*}
	where in the last inequality we have used~\cite[Lem.~3.3]{high:ASNA2},
	and so the proof is completed by induction. 
\end{lemma}

Based upon Lemma~\ref{thm:xpower-err}, we obtain the following result, which bounds the forward error of a polynomial formed by assembling matrix powers that have already been computed via repeated multiplication. 
It can be viewed as a specification of the more general result~\cite[Thm.~4.5]{high:FM}.
The proof is similar to that of Theorem~\ref{thm:ps-mp-err-norm} by induction and is thus omitted.

\begin{theorem}\label{thm:xpoly-err}
If the first $t$ positive powers 
of $X$ are formed by repeated multiplication in precision $u$ (with $\wh{X}_t = \fl(X^t)$ denoting the computed matrix) and
$\psi=\sum_{j=0}^{t}a_{j}X^j$ is evaluated in precision $u$ by
\begin{code}
        $\widehat{\varphi}_0=\fl(a_{0}I)$ \\
        for $j=1:t$ \\
        \> $\widehat{\varphi}_j=\fl(\widehat{\varphi}_{j-1} + \fl(a_{j}\widehat{X}_{j}))$ \\
        end \\
        $\widehat{\psi}=\widehat{\varphi}_t$  	
\end{code}
then the computed $\wh{\psi}$ satisfies
\begin{equation}\label{eq:xpoly-err}
        |\wh{\psi}-\psi(X)| \le 
        \gamma_{t}|a_{0}| I + \sum_{j=1}^{t}\gamma_{(j-1)(n-1)+t+1} |a_{j}| 
        |X|^j.
\end{equation}	
\end{theorem}
Since potentially arbitrarily weak inequalities such as $|X^i|\le|X|^i$ are used in the derivation, the bound~\eqref{eq:xpoly-err} can be unduly  pessimistic, but as an a priori 
bound it cannot be improved without further assumptions.
The bound is immediately applicable to the computed polynomial
$\wh{B}_{\nu-1}$ in precision $u=u_{\nu-1}$, and we have 
\begin{align}\label{eq:B0-err-bnd}
	|\wh{B}_{\nu-1}-B_{\nu-1}(X)| &\le 
	\gamma_{s-1}|b_{s(\nu-1)}| I + \sum_{j=1}^{s-1}\gamma_{(j-1)(n-1)+s} |b_{s(\nu-1)+j}| 
	|X|^j \nonumber \\
	&\le \gamma_{(s-2)(n-1)+s} \sum_{j=0}^{s-1} |b_{s(\nu-1)+j}| |X|^j 
	=: \gamma_{(s-2)n+2}\widetilde{B}_{\nu-1}(|X|),
\end{align}	
where $\widetilde{B}_{\nu-1}(X)=\sum_{j=0}^{s-1} \abs{b_{s(\nu-1)+j}} X^j$. Thus, generally, a sufficient condition for~\eqref{eq:cond-B0} to hold is  
$\norm{B_{\nu-1}(X)}\approx \norm{\widetilde{B}_{\nu-1}(|X|)}$, which is true if there is no significant cancellation in forming $B_{\nu-1}(X)$. This is the case, for example, when the coefficients $b_{s(\nu-1)+j}$, $j=0\colon s-1$ have the same sign and $X>0$, or, when the $\abs{b_{s(\nu-1)+j}}$ decay rapidly and there is no significant cancellation in forming the first few terms of $B_{\nu-1}(X)$.

On the other hand, since all the required powers of $X$ are formed in precision $u$, 
we have $\norm{\widehat{Y} - Y}\le \gamma_{(s-1)n}\norm{X}^{s}$ from Lemma~\ref{thm:xpower-err}. 
Therefore, 
from~\eqref{eq:ui-ui-1-relation} we have
\begin{equation}\label{eq:cond-Y}
	\norm{\widehat{Y} - Y}\lesssim snu \norm{X}^s = sn\tau_\nu u_\nu \norm{X}^s \lesssim u_\nu\norm{X^s}
\end{equation}
if $sn\tau_\nu\norm{X}^s \lesssim \norm{X^s}$. 
In any case, the validity of this relation will depend on the matrix $X$ and $\tau_\nu$. 
A special instance is when $X\ne 0$ is nilpotent with index $s$ (so $X^s=0$), where the condition 
$\norm{\wh{Y}-Y}\le u_\nu\norm{X^s}$
cannot possibly be satisfied ($\widehat{Y}$ can contain significant rounding errors), but we are not interested in this case where the evaluation of~\eqref{eq:ps-Horner} becomes trivial because $p_m(X)=B_0$.

\section{Taylor approximants to the matrix exponential}\label{sect:exp}
In this section we consider the concrete setting when the matrix polynomial $p_m(X)$ of~\eqref{eq:ps} is the truncated Taylor approximant of order $m$ to the
matrix exponential of $X$, where the coefficients $b_i=1/i!$ decay
super-exponentially. 

We can show that if the $1$-norm of $X$ satisfies 
\begin{equation}\label{eq:X-1norm-bound}
	\normi{X}\le s/\eu, \quad 1\le s\le m,
\end{equation}
where $\eu\approx 2.718$ is Euler's constant and $s\in\mathbb{N^+}$ is the parameter in the PS scheme,
then in general the conditions~\eqref{eq:Bi-Xs-bnd} are satisfied with $\nu=1$ and accuracy of the computed polynomial $\widehat{p}_m$ of $p_m(X)$ is guaranteed. The condition~\eqref{eq:X-1norm-bound} can be too stringent on the input matrix to the polynomial, which, when the polynomial is embedded into scaling and squaring method for the matrix exponential~\cite[sect.~10.3]{high:FM}, can cause overscaling issue and lead to loss of accuracy in the squaring steps,
but the analysis can give us some insight into the feasibility of the conditions~\eqref{eq:Bi-Xs-bnd} and the behaviour of our mixed-precision Paterson--Stockmeyer algorithm in floating point arithmetic.

We will later discuss why we think the constraint~\eqref{eq:X-1norm-bound} can be relaxed in practice: disregarding the constraint does not prevent the conditions~\eqref{eq:Bi-Xs-bnd} being satisfied with some $\nu\in\left[1,r\right]$ (often $\nu=1$) in the practically used algorithms; and, on the other hand, from our experience and the experiments in the existing literature, it does not appear to be
harmful for the accuracy of a practical algorithm in the occurrence of rounding errors.
The caveat is that there is no theoretical guarantee of the accuracy in the part of the PS scheme done in the working precision
as our analysis only accounts for the potential use of lower precisions.

\subsection{Analysis with the norm constraint}
Suppose $\normi{X}=\sigma$ for some $0<\sigma\le s/\eu$. 
We noticed $\sqrt[s]{s!}$ can be very well approximated by $s/\eu+1$ for $s$ in a practical range, say, $s\le 25$. (Note that for large $s$ Stirling's approximation gives $\sqrt[s]{s!}\sim s\sqrt[2s]{2\pi s}/\eu \sim s/\eu$.) 
Therefore, we have 
\begin{equation}\label{eq:Xs-1norm-bnd}
	\normi{X^s}\le \normi{X}^s= \sigma^s \le (s/\eu)^s\le s!.
\end{equation}
It follows that, for $i=2\colon r$,
\begin{align*}
	\tau_i=\frac{\normi{B_{i}}\normi{X^s}}{\normi{B_{i-1}}} &=
	\frac{\left\Vert\frac{1}{(is)!}I+\frac{1}{(is+1)!}X+\cdots+\frac{1}{(is+s-1)!}X^{s-1}\right\Vert_1\left\Vert X^s\right\Vert_1}
	{\left\Vert\frac{1}{((i-1)s)!}I+\frac{1}{((i-1)s+1)!}X+\cdots+\frac{1}{((i-1)s+s-1)!}X^{s-1}\right\Vert_1} \\
	&\le
	\frac{\left(\frac{1}{(is)!}+\frac{\sigma}{(is+1)!}+\cdots+\frac{\sigma^{s-1}}{(is+s-1)!}\right)s!} 
	{\frac{1}{((i-1)s)!}-\left(\frac{\sigma}{((i-1)s+1)!}+\frac{\sigma^2}{((i-1)s+2)!}+\cdots+\frac{\sigma^{s-1}}{((i-1)s+s-1)!}\right)}\\
	&\le
	\frac{\frac{1}{(is)!}\left(1+\frac{\sigma}{is+1}+\cdots+\frac{\sigma^{s-1}}{(is+1)^{s-1}}\right)s!}
	{\frac{1}{((i-1)s)!}-\frac{\sigma}{((i-1)s+1)!}\left(1+\frac{\sigma}{(i-1)s+2}+\cdots+\frac{\sigma^{s-2}}{((i-1)s+2)^{s-2}}\right)}=:\gamma(s,i).
\end{align*}
Since $\sigma\le s/\eu$ implies, for $i\ge 2$, 
$$ 
r_1:=\frac{\sigma}{(i-1)s+2}< \frac{1}{\eu}<1,\quad r_2:=\frac{\sigma}{is+1}< \frac{1}{2\eu}<1,
$$
it follows that, for $r_1$ and $r_2$ sufficiently close to zero (which is true when $\sigma\approx 0$ or when $s$ or $i$ is large),
\begin{align}\label{eq:gamma-s}
	\gamma(s,i)&=
	\frac{\frac{s!}{(is)!}\cdot\frac{1-r_2^s}{1-r_2}}
	{\frac{1}{((i-1)s)!}-\frac{\sigma}{((i-1)s+1)!}\cdot\frac{1-r_1^{s-1}}{1-r_1}} \approx
	\frac{s!(is-s)!}{(is)!\left(1-\frac{\sigma}{(i-1)s+1}\right)} \nonumber\\
	&\le  \left(1-\frac{1}{\eu(i-1)}\right)^{-1}{is\choose s}^{-1} \le \frac{\eu}{\eu-1}\cdot\frac{s^s}{(is)^s}	= \frac{\eu}{\eu-1}i^{-s}
	\approx 1.58i^{-s}.
\end{align}
This shows that, for a chosen $s$, $\tau_i\le\gamma(s,i)$ (recall that $\tau_i$ is from~\eqref{eq:Bi-Xs-bnd}) decreases at least polynomially as $i$ increases, and that, for a fixed $i\ge 2$, $\tau_i$ decreases at least exponentially as $s$ increases.

We have 
\begin{equation}\label{eq.bnd-B1}
	\normi{B_1} 
	\le  \frac{1}{s!}\biggl(1+\frac{\sigma}{s}+\cdots+\frac{\sigma^{s-1}}{s^{s-1}}\biggr)
	= \frac{1-(\sigma/s)^s}{s!(1-\sigma/s)}\le \frac{1}{(1-\sigma/s)s!}.
\end{equation}
Both the bound and $\normi{B_1}$ tend to $1/s!$ as $\sigma\to 0$, which shows this bound tends to be a good approximation to $\normi{B_1}$ as $\normi{X}$ tends to zero.
Bound~\eqref{eq.bnd-B1} implies $\normi{B_1}\le \frac{1}{(1-1/\eu)s!}$, and by using this relation and the approximation $(s/\eu)^s\le s!$, we arrive at
\begin{align}\label{eq:B0-B1-tau-cond}
	\tau_1 =\frac{\normi{B_1}\normi{X^s}}{\normi{B_0}} &\le \frac{\eu\normi{X}^s}{(\eu-1)s!\normi{B_0}}\cdot\frac{\normi{X^s}}{\normi{X}^s} \nonumber \\
	&\lesssim
	\frac{(s/\eu)^s}{s!\normi{B_0}}\cdot\frac{\normi{X^s}}{\normi{X}^s} \le
	\frac{1}{\normi{B_0}}\cdot \frac{\normi{X^s}}{\normi{X}^s},
\end{align}
where $\normi{B_0}\approx\normi{\eu^X}\ge \eu^{-\normi{X}}$~\cite[Thm.~10.10]{high:FM} and 
$\normi{X^s}/\normi{X}^s$ is bounded above by $1$ but can be arbitrarily small. 
This shows that $\tau_1$ is bounded above by a quantity which tends to zero as $s$ increases.

Now consider the effects of rounding errors on the evaluation of $Y=X^s$ and
$B_0$ in precision $u=u_0$ (as now we consider $\nu=1$).
From~\eqref{eq:cond-Y} we need to check whether we have $sn\tau_1\norm{X}^s \lesssim \norm{X^s}$. We have, using the approximation $B_0(X) \approx \eu^{X}$ from~\eqref{eq:Bi-Xs-bnd} and the bounds~\eqref{eq:Xs-1norm-bnd} and~\eqref{eq.bnd-B1},
\begin{equation*}
	\frac{sn\tau_1\normi{X}^s}{\normi{X^s}} = \frac{sn\normi{B_1}\normi{X}^s}{\normi{B_{0}}} \\
	\lesssim \frac{sn\cdot s!}{\normi{\eu^{X}}}\cdot 
	\frac{1}{(1-\sigma/s)s!} \\	
	\lesssim   sn\eu^{\normi{X}},
\end{equation*}
which shows, given that $\normi{X}$ is nicely bounded, $sn\tau_1\normi{X}^s$ is approximately bounded above by a mild multiple of $\normi{X^s}$, and therefore, we should expect $Y=X^s$ to be evaluated to satisfying accuracy.
On the other hand, since $B_0$ has all positive coefficients, it follows from~\eqref{eq:B0-err-bnd} that
\begin{equation*}
	\normi{\wh{B}_0-B_0(X)} \le \gamma_{(s-2)n+2}B_0(\normi{X}) 	\approx  \gamma_{(s-2)n+2}\eu^{\normi{X}},
\end{equation*}
and we deduce, using~\cite[Thm.~10.10]{high:FM}, 
\begin{align*}
	\normi{\wh{B}_0-B_0(X)} &\lesssim \gamma_{(s-2)n+2}\eu^{\normi{X}} \cdot \eu^{\normi{X}}\eu^{-\normi{X}} \\
	&\le \gamma_{(s-2)n+2}\eu^{\normi{X}}\cdot  \eu^{\normi{X}}
	\normi{\eu^X} \\	
	&\approx  \gamma_{(s-2)n+2}\normi{B_0(X)}\eu^{2\normi{X}} 
	\le \gamma_{(s-2)n+2}\normi{B_0(X)}\eu^{2s/\eu}.
\end{align*}
Hence the relative error in $\wh{B}_0$ is bounded approximately by $\gamma_{(s-2)n+2}\eu^{2s/\eu}$, which is a satisfactory bound for practical values of $s$.
We have empirically found that $\wh{B}_0$
is typically computed to close to full working precision (the relative error in
$\wh{B}_0$ is typically close to $u$) for matrices of varying size generated pseudo-randomly and from the MATLAB gallery. This is consistent with the analysis which shows the rounding errors in the evaluation of $B_0$ are nicely bounded, and it is also a possible consequence of the fact that 
the underlying basic linear algebra subprograms (BLAS) in MATLAB 
uses blocked algorithms to reduce the error constant~\cite{high21n}.

\subsection{Applicability in the scaling and squaring algorithms for the matrix exponential}
The discussion in the previous subsection implies that we could build a mixed-precision PS algorithm for the $m$th-order Taylor approximant to the matrix exponential under the constraint $\normi{X}\le s/\eu$, where the parameter $s$ could exceed $\ceil{\sqrt{m}}$ in order for the desired decaying property of the polynomial coefficients~\eqref{eq:Bi-Xs-bnd} to hold, and accuracy of the algorithm could in general be guaranteed.
However, in the stat-of-the-art algorithms for the matrix exponential~\cite{alhi09a},~\cite{fahi19}, which employ the scaling and squaring idea
\begin{equation}\label{eq:scal-idea}
	\eu^A=\left(\eu^{X}\right)^{2^{\ell}}\approx p_m(X)^{2^{\ell}},\quad \ell\in\mathbb{N^+}\cup\{0\},
\end{equation}
the admittable (scaled) matrix $X:=2^{-\ell}A$ can have a $1$-norm that does not satisfy the constraint $\normi{X}\le s/\eu$. 
This is because the thresholds for accepting certain $X$ in these algorithms are determined by forward or backward truncation error bounds of the approximant to the exponential on the scaled matrix in exact arithmetic, and these thresholds in fact disregard the rounding errors in the computation of the approximant $p_m(X)$.
For example, the Taylor-based algorithm of~\cite{fahi19} 
requires that $X$ satisfy
\begin{equation}\label{eq:alpha_m-X}
	\normi{\eu^X-p_m(X)}\le \abs{\eu^{\alpha_m(X)}-p_m(\alpha_m(X))}\le u\xi(X),
\end{equation}
where $\xi(X)$ is some practical estimate to $\normi{\eu^X}$ and 
\begin{multline*}
	\alpha_m(X) =  \max\bigl( 
	\normi{X^{d^*}}^{1/d^*},\normi{X^{d^*+1}}^{1/(d^*+1)} \bigr), \\
		d^* = \max_d \{d\in\mathbb{N}^+\colon d(d-1)\le m+1 \} = 
	\left\lfloor \frac{1 + \sqrt{4m + 5}}{2}  \right\rfloor.
\end{multline*}

In principle, the constraint $\normi{X}\le s/\eu$ does not prevent the potential algorithm from being
embedded into any existing scaling and squaring algorithm based on the Taylor approximants, for example, those employ the $\alpha_m$-based bound~\eqref{eq:alpha_m-X}. This is because, for a scaled matrix $X$ accepted by one of such algorithms with $\normi{X}>s/\eu$, one could always further scale $X$ to $Z=2^{-\ell_0}X$, $\ell_0\in\mathbb{N}^+$ such that $\normi{Z}\le s/\eu$ and $Z$ remain admittable by the algorithm since $\alpha_m(Z)=2^{-\ell_0}\alpha_m(X)<\alpha_m(X)$. 
In this way, instead of invoking~\eqref{eq:scal-idea}, the algorithm would be using the approximation
$$
\eu^A=\left(\eu^{2^{-(\ell+\ell_0)}A}\right)^{2^{\ell+\ell_0}}\approx p_m(Z)^{2^{\ell+\ell_0}},
$$ 
which, from our discussion in the previous subsection, has more refined bound on the rounding errors in 
the computed approximant $\widehat{p}_m$. The algorithm nevertheless can require substantially more squaring steps in the final squaring phase, which is very sensitive to rounding errors~\cite[p.~247]{high:FM}, because 
$\alpha_m(X)$ can be much smaller than $\normi{X}$ for nonnormal $X$~\cite{alhi09a} and 
a matrix $X$ accepted by an $\alpha_m$-based bound can have huge $1$-norm. For example, consider the matrix 
\begin{equation}\label{eq:nonnormal-mat-A}
	A = \begin{bmatrix}
		-0.1 & 10^6 \\
		0  & -0.1
	\end{bmatrix},
\end{equation}
for which the $1$-norms of the powers of $A$ decay exponentially and the Taylor-based algorithm of~\cite{fahi19} with $u=10^{-64}$ chooses $m=42$ and accepts $X=A/2$ ($\ell=1$), despite the large $(1,2)$ element.
In this case $\alpha_m(X)\approx 0.66$ but $\normi{X}=5\times10^5$ and with $s=\ceil{\sqrt{m}}$ the number of extra squarings required is
$
\ell_0 \ge \left\lceil{\log_2\left(\eu\normi{X}/s\right)}\right\rceil = 18.
$

In fact, the scaling and squaring algorithms for the matrix exponential~\cite{alhi09a}, \cite{bbc19}, \cite{fahi19}, \cite{sidr15},
which disregard the occurrence of rounding errors in the computed approximant $\widehat{p}_m$ when determining the thresholds for accepting the scaled matrix,
have been observed to work well in practice, even for $X$ with a large $1$-norm but a small or moderate $\alpha_m(X)$ associated with the chosen $m$.
Disregarding the constraint $\norm{X}\le s/\eu$, which can cause overscaling issue for the scaling and squaring algorithms, does not appear to be
harmful for the accuracy of a practical algorithm in the occurrence of rounding errors.
Then in this case the question is to what extent~\eqref{eq:Bi-Xs-bnd} can be satisfied? We have found the condition often holds with $\nu=1$ if $\norm{X^s}$ is small, which is consistent with the discussion following~\eqref{eq:Bi-Xs-bnd}.
In the scaling and squaring algorithms that employ the $\alpha_m$-based bound~\eqref{eq:alpha_m-X}, for example,~\cite{fahi19}, the $\alpha_m(X)$ can be very small (associated with the chosen $m$) on some tested matrices even if
$\norm{X}$ is large, in which case the value of $\normi{X^{d^*}}$ is necessarily small, where $d^*\approx \sqrt{m}$ just matches the default parameter $s=\floor{\sqrt{m}}$ or 
$s=\ceil{\sqrt{m}}$ in the fix-precision PS scheme, so~\eqref{eq:Bi-Xs-bnd} often holds with $\nu=1$.
In the least preferred case where~\eqref{eq:Bi-Xs-bnd} is not met for any $\nu\in\left[1,r\right]$, still, we can simply compute $q(X)$ from~\eqref{eq:ps-Horner-part} in the working precision, in which case the algorithm recovers the fix-precision PS scheme.
\begin{codefragment}[t]
	\caption{Computing $B_0$ and $Y=X^s$ in precision $u$.}
	\label{alg.exp.B0}
	\Fn{\textsc{Eval$B_0$}$(X$, $s\in\mathbb{N}^+$$)$\\
		\funcomment{Form the first $s$ positive powers of $X$ in $\X$ and then
			compute $B_0 = \sum_{j=0}^{s-1}X^j/j!$ and $Y=X^s$ using
			elements of $\X$.}}{
		$B_0 \gets I$, $\X_1 \gets X$\;
		\For{$j \gets 1$ \algto $s-1$}{
			$B_0 \gets B_0 + \X_{j}/j!$ \;
			$\X_{j+1} \gets \X_{j} X$ \;
		}
		$Y \gets \X_{s} $  \;
		\Return{$\X$, $B_0$, $Y$}
}\end{codefragment}

\begin{algorithm2e}
	\setstretch{1}
	\caption{Mixed-precision Paterson--Stockmeyer scheme for the Taylor approximants of the matrix exponential.}
	\label{alg.mpps}
	\SetAlgoVlined
	{\nonl
		\begin{minipage}{0.95\linewidth}
			\pretolerance=100\tolerance=200\hyphenpenalty=5
			Given $X\in\mathbb{C}^{\nbyn}$ this algorithm computes an $m$th order Taylor approximant 
			$P\equiv p_m(X)$ in the form of~\eqref{eq:ps} using the Paterson--Stockmeyer scheme in floating-point arithmetic. The algorithm starts with the user-specified precision $u$ and potentially uses multiple lower precisions $u_i\ge u$
			aiming to produce a relative error of order $nu$.
			The logical parameter \texttt{fix\_params} determines whether some $s>\ceil{\sqrt{m}}$ is allowed to be used for small-normed matrices.
			The pseudocode of \textsc{Eval$B_0$} is given in Fragment~\ref{alg.exp.B0}.
	\end{minipage}}
	\BlankLine
	$s \gets \ceil{\sqrt{m}}$ \; 
	$u_0 \gets u$ \;
	$[\X, B_0, Y] \gets \textsc{Eval$B_0$}(X, s)$ \;
	\uIf{$\texttt{fix\_params}$}{\label{alg.mpps.line.fix.params-start}
		$r \gets \floor{m/s}$ \;
		$\nu \gets r+1$ \;
		\For{$i \gets 1$ \algto $r$}{
			Assemble $B_i$ using elements in $\X\cup\{I\}$ and then estimate $\normi{B_i}$\;
			$u_i \gets \normi{B_{i-1}}u_{i-1}/(\normi{B_i}\normi{Y}^i)$\;
		}
		$\nu\gets \min\{i\colon u_i\ge \delta u\}$, $u_{\nu-1},u_{\nu-2},\dots,u_{1}\gets u$ \;
		}\label{alg.mpps.line.fix.params-end}
	\uElseIf{$\normi{X}\le s/\eu$}{\label{alg.mpps.line.var.params-start}
		\While{$\normi{Y}> \normi{B_0}\normi{X}^s$ \algand $s< m$}{
			$B_0 \gets B_0+Y/s!$ \;
			$s\gets s+1$ \;
			$\X_s\gets XY$ in precision $u_0$ \;
			$Y\gets \X_s$ \;}
		$r \gets \floor{m/s}$ \;
		\For{$i \gets 1$ \algto $r$}{
			Assemble $B_{i}$ using elements in $\X\cup\{I\}$ and then estimate $\normi{B_{i}}$\;
			$u_{i} \gets \normi{B_0}u_0/(\normi{B_{i}}\normi{Y}^{i})$\;
		}
		\label{alg.mpps.line.var.params-end}
	}
	\Else{
		$\textit{\texttt{fix\_params}} = \texttt{true}$ and go to line~\ref{alg.mpps.line.fix.params-start}\;
	}
	$P = B_r$ \;
	\For{$i \gets r$ \algdownto $1$}{
		Compute $P\gets PY$ in precision $u_i$\;
		Form $P\gets P + B_{i-1}$ in precision $u_{i-1}$ \;	}
	\Return{$P$}\;
	\SetAlgoNoLine
\end{algorithm2e}

\subsection{The mixed-precision Paterson--Stockmeyer algorithm}
\label{sect.mpps.exp.taylor}
Summarising the discussion in the previous sections, we can build a mixed-precision PS algorithm which starts with computing $B_0$ straightly in precision $u$ with the default parameter $s=\ceil{\sqrt{m}}$ and then proceeds differently depending on $\normi{X}$.

If $\normi{X}\le s/\eu$, the algorithm can be made to
find a sufficiently large value of $s$ such that lower-than-the-working precisions can be used in the matrix products associated with $B_{i}$, $i=1\colon r$. It
increments $s$ and updates $B_0$ and $Y$ until the bound $\tau_1\le \normi{X^s}/\left(\normi{B_0}\normi{X}^s\right)\le1$ from~\eqref{eq:B0-B1-tau-cond} is satisfied. 
Here one can check~\eqref{eq:Bi-Xs-bnd} only for $i=1$ since
we have shown that the $\tau_i$, $i=2\colon r$ tend to decay at least polynomially as $i$ increases (see~\eqref{eq:gamma-s}) and we found if the first condition of~\eqref{eq:Bi-Xs-bnd} is satisfied, then the remaining conditions therein are met for $\normi{X}\le s/\eu$. 

If $\normi{X}>s/\eu$, there is no theoretical support that condition~\eqref{eq:Bi-Xs-bnd} will hold with $\nu=1$ for some $s\le m$, meaning that  
we might not be able to safely use lower-than-the-working precisions in forming
the matrix products associated with $B_{1},B_{2},\dots$ for any feasible choice of $s$.
The algorithm therefore sticks with $s=\ceil{\sqrt{m}}$ to form the $B_i$, $i=1\colon r$, and calculates the precision $u_i$ according to the norms of the coefficient matrices: it sets $u_i$ by $\max\left(u,\mbox{right-hand side of}~\eqref{eq:choose-ui}\right)$  
until $u_i\ge \delta u$ (say, $i=\nu$), so lower precisions $u_r\ge u_{r-1}\ge\cdots\ge u_{\nu}$
are exploited in the computation of $q(X)$ and the remaining part of $p_m(X)$ is then computed via~\eqref{eq:pm-prec-u-part} in the working precision $u$.
By default we set $\delta=10$, so the algorithm will switch to a lower precision when appropriate, even if the number of significant digits decreases by just $1$.

After $B_i$ and $u_i$ are computed for all $i$ the algorithm then executes the Horner's method~\eqref{eq:ps-Horner} for $p_m(X)$ with the matrix products and sums done in the appropriate precisions. 
The complete mixed-precision PS algorithm is presented in Algorithm~\ref{alg.mpps}. The algorithm takes the matrix $X\in\mathbb{C}^{\nbyn}$, the order $m$ of the used Taylor approximant, and the working precision $u$ as input arguments, where we have also introduced 
a logical parameter \texttt{fix\_params} that determines whether some $s>\ceil{\sqrt{m}}$ is allowed to be used for small-normed matrices,
for which there is in principle a balance: fixing the parameter $s=\ceil{\sqrt{m}}$ minimizes the total number of matrix multiplications done in all precisions but the algorithm might need perform some matrix multiplications associated with $B_1, B_2, \dots$ in the working precision $u$, while using a larger $s>\ceil{\sqrt{m}}$ requires some extra number of matrix products in total yet it might enable the 
low-than-the-working precisions $u_1,u_2,\dots,u_r$ chosen by the algorithm to be even lower. 
Note the Boolean value of \texttt{fix\_params} does not affect the executed steps of the algorithm for input matrices with $\normi{X}> s/\eu$ as the default $s=\ceil{\sqrt{m}}$ is used in this case.

The economics of the mixed-precision PS scheme are different from Algorithm~\ref{alg.fpps}, the fixed-precision PS: computing the matrix products is still computationally the dominant part, but there is a discrimination between these products as they are performed in different precisions and those done in the working precision are the most expensive.
Overall, if lines~\ref{alg.mpps.line.var.params-start}--\ref{alg.mpps.line.var.params-end} are invoked in Algorithm~\ref{alg.mpps}, then it requires
$\ceil{\sqrt{m}}-1\le s-1\le m-1$ matrix multiplications in precision $u$ and one matrix multiplication in each of $u_i>u$, $i=1\colon r$, where
$1\le r = \floor{m/s}\le\ceil{\sqrt{m}}$ (when $s=m$ the PS scheme actually
degenerates to evaluation via explicit powers and hence no matrix
multiplications are formed in $u_i$).
Otherwise, when
lines~\ref{alg.mpps.line.fix.params-start}--\ref{alg.mpps.line.fix.params-end} are executed, the dominant cost of the algorithm is 
$\ceil{\sqrt{m}}+\nu-2$ matrix multiplications in precision $u$
and one matrix multiplication in each of $u_{\nu},u_{\nu+1},\dots,u_{r}$, 
where $1\le \nu\le r$ and a smaller $\nu$ implies more matrix products are performed in precisions lower than $u$; and in the case where
$u_i\ge \delta u$ does not hold for any $i$, 
Algorithm~\ref{alg.mpps} recovers the 
fix-precision PS scheme (with some additional norm estimations), which is clearly the worst case for the algorithm in terms of efficiency.
On the other hand, in the case of optimal
efficiency, Algorithm~\ref{alg.mpps} requires 
$\ceil{\sqrt{m}}-1$ matrix multiplications in $u$, which is only approximately half of the
matrix multiplications required by the fix-precision PS scheme.

There is also the computational overhead from $O(n^2)$ flops in Algorithm~\ref{alg.mpps} resulting from the norm computations and it can be nonnegligible for small size problems.
However, since the norm estimations in the algorithm are needed only for determining the precisions $u_i$, a rough estimation to correct order of magnitude is sufficient and this can be done in a low precision, with the mitigation that the precision of norm estimation can be much lower than the working precision as the algorithm mainly targets arbitrary precision environment.

Finally, we remark that Algorithm~\ref{alg.mpps} is readily employable by the Taylor-based scaling and squaring algorithms for the matrix exponential, and its framework can be adjusted with little modification for the computation of other matrix polynomials with scalar coefficients that quickly decay in modulus, such as the polynomials in the numerator and denominator of the Pad\'e approximants of exponential-like functions.
We will examine the generality of the framework in Section~\ref{sect.general.framework}.

\subsection{Numerical experiments}\label{sect.numerical.exp.taylor}

 
All our experiments were run using the 64-bit GNU/Linux version of MATLAB 23.2 (R2023b) on a desktop computer
equipped with an Intel i5-12600K processor running at 3.70 GHz
and with 32GiB of RAM;
Two different arithmetics are tested in the experiments. The main one is \textit{arbitrary precision arithmetic} simulated by the Advanpix Multiprecision Computing Toolbox (Version~5.1.1.15444)~\cite{adva-mct}, where the algorithms start with the working precision $u$ and can potentially use internally multiple lower precisions $u_i\ge u$, $i=1\colon r$ that can be arbitrarily chosen. We also tested the setting of arithmetic where only IEEE double, IEEE single, and bfloat16 half precision (by Google Brain\footnote{\url{https://research.google}})
are available, the last being simulated by the \t{chop}\footnote{\url{https://github.com/higham/chop}} function~\cite{hipr21}.
The MATLAB code that produces the test results is available from GitHub.\footnote{\url{https://github.com/xiaobo-liu/mp-ps}}

As aforementioned, the mixed-precision PS method Algorithm~\ref{alg.mpps} mainly aims for computing the matrix polynomials required by the Taylor-based scaling and squaring algorithm for the matrix exponential, in which case the argument matrix $X$ is often pre-scaled by $X=2^{-\ell}A$ (see~\eqref{eq:scal-idea}), with the scaling parameter $\ell$ (in fact, together with the polynomial degree $m$) set by the scaling and squaring algorithm. 
Therefore, to simulate a practical experimental setting, in the experiments within this section we feed the matrix $A$, of which the exponential is supposedly of interest, and the working precision $u$ to 
the state-of-the-art Taylor-based scaling and squaring algorithm aiming for arbitrary precision environment~\cite[Alg.~4.1]{fahi19} to get the test data $X$ and the polynomial degree $m$.
We then compare Algorithm~\ref{alg.mpps} with its \textit{fixed-precision counterpart} (Algorithm~\ref{alg.fpps} specified to the Taylor approximant of the matrix exponential) in precision $u$ for computing the matrix polynomial $p_m(X)$, the latter being exactly the PS algorithm employed in~\cite{fahi19}.
Note that, when we talk about any test matrices in the following experiments, we are referring to $A$, the matrix prior to the scaling.

In all experiments the reference solution is computed by running Algorithm~\ref{alg.fpps} (specified to the corresponding approximant in the test) using the Advanpix Toolbox with twice the decimal digits of the working precision. The logical parameter \texttt{fix\_params} is switched on (for small-normed matrices) in
Algorithm~\ref{alg.mpps} since we found virtually no difference in terms of accuracy between the two variants with \texttt{fix\_params} switched on or off, while the one with $\texttt{fix\_params}=\texttt{true}$ is usually faster. This is possibly a sign that the benefits brought by minimizing the total number of matrix products outweighs performing the lower-precision matrix products using fewer digits in our implementation.

\subsubsection{Distributions of the used lower-than-working-precisions}

\begin{table}[t]
	\centering\footnotesize
	\caption{The approximant degree $m$ is chosen automatically by 
	 the Taylor-based algorithm of~\cite{fahi19}
	 given the input matrix $X$ and the precision $u$. The equivalent decimal digits of the working precisions are reported in the first column. The $d_i$ represents the equivalent decimal digits of precision $u_i$ and $C_r$ is the approximate cost ratio of Algorithm~\ref{alg.mpps} over its fixed-precision counterpart in precision $u$.}
	\input{matlab/tabs/table_cauchy_100.tex}
	\label{tab.exp_cauchy}
\end{table}

First, to get some insight into the distribution of the lower-than-the-working precisions chosen by Algorithm~\ref{alg.mpps}, 
we take $A=\ $\verb|gallery('cauchy',n)| with $n=100$ and 
compute the Taylor approximant from its matrix exponential with several choices of the working precision $u$;
for this matrix it has $\normi{A}\approx 4.20$ and the Taylor-based algorithm of~\cite{fahi19} chooses the scaling parameter $\ell=0$ (no scaling and thus $X=A$). 
Some important algorithmic characteristics on the matrix are reported in
Table~\ref{tab.exp_cauchy}.

We see from the table that 
$\tau_i=u_{i-1}/u_i$ is in general decreasing, which is consistent with our analysis (see~\eqref{eq:gamma-s}). 
Recall that the algorithm requires totally $s-1+r$ matrix multiplications with the default $s=\ceil{\sqrt{m}}$, and we see that typically about a fifth of the matrix multiplications were performed in precision $u^{1/2}$ or much lower in all tested cases. 
If the algorithmic complexity (measured in number of matrix multiplications) is assumed to be linearly proportional to the number of digits used, then we can calculate 
the approximate cost ratio of Algorithm~\ref{alg.mpps} over its fixed-precision counterpart in precision $u$ as
\begin{equation}\label{eq:Cr}
	C_r = \frac{(s-1)\log_{10}u + \sum_{i=1}^{r}d_i}{(s+r-1)\log_{10}u}.
\end{equation}
Note that the assumption on the algorithmic complexity is realistic yet slightly pessimistic because scalar multiplications and divisions scale even faster with the number of digits.
The ratio of cost $C_r$ ignores the operations that depend quadratically or linearly on $n$ and 
does not consider other important factors for the performance, e.g., data transfer, memory access, and precision switching;
in this sense the measure is quite theoretical and in general should not be  interpreted as an implication for the reduction of the actual computational time.

\begin{figure}
%
		\begin{subfigure}{1\linewidth}
			\centering
			\includegraphics[height=4.3cm]{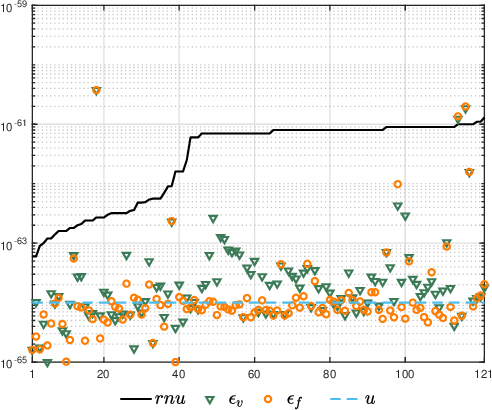} 
			\includegraphics[height=4.3cm]{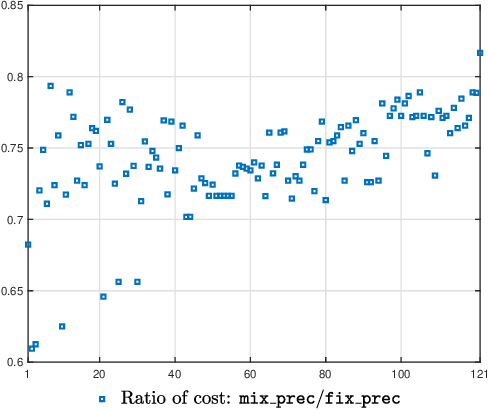}
			\caption{$u = 10^{-64}$.}
		\end{subfigure} \\
				\vspace{-15pt}
				
		\begin{subfigure}{1\linewidth}
			\centering
			\includegraphics[height=4.3cm]{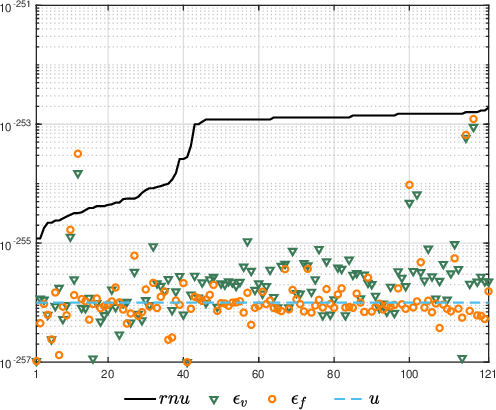}
			\includegraphics[height=4.3cm]{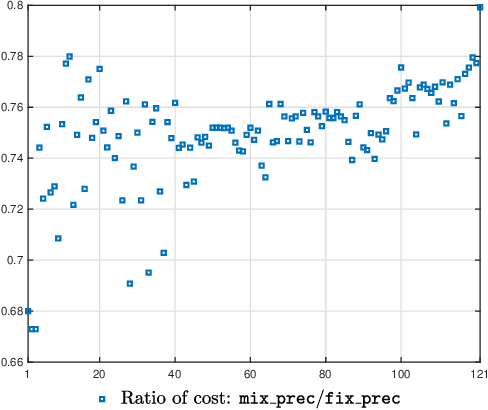}
			\caption{$u = 10^{-256}$.}
		\end{subfigure}
		\vspace{-15pt}
		
		\caption{Left: The relative error $\epsilon_v=\normi{\widehat{p}_m-{p}_m(X)}/\normi{{p}_m(X)}$ 
			produced by Algorithm~\ref{alg.mpps} compared with the relative error $\epsilon_f$ produced by the fixed-precision counterpart from~\cite{fahi19} with $s=\ceil{\sqrt{m}}$ on various matrices with $2\le n\le 100$. Right: The associated approximate ratio of cost $C_r$ in~\eqref{eq:Cr}.}
		\label{fig.exp_taylor_ap}
\end{figure}

\subsubsection{Accuracy and theoretical speed-up}\label{sect.numerical.accuracy.speedup}

Next, we compare Algorithm~\ref{alg.mpps} with its fixed-precision counterpart, using the test matrices 
from the Anymatrix collection~\cite{himi-am},~\cite{himi21} and from various literature of matrix exponentials that are used in~\cite{alhi09a},~\cite{fahi19};  
matrices from the latter are collected in an Anymatrix group accessible on GitHub.\footnote{\url{https://github.com/xiaobo-liu/matrices-expm}}
There are in total $121$ test matrices of size ranging from $2$ to $100$.

We see from Figure~\ref{fig.exp_taylor_ap}~$(a)$--$(b)$ that the relative errors produced by Algorithm~\ref{alg.mpps} are comparable to that of the fixed-precision counterpart and no larger than $rnu$
in most cases. On the other hand, we see that in almost all cases the computational cost of the algorithm is in theory at least 20\% lower compared with its fixed-precision counterpart on the test set, and in few cases the savings approach 40\%.

\subsubsection{Code profiling and runtime comparison}

\begin{table}
	\renewcommand{\tabcolsep}{4pt}
	\footnotesize
	\caption{Execution time breakdown of Algorithm~\ref{alg.mpps} and its fixed-precision counterpart, run in 64 and 256 decimal significant digits on three classes of matrices of increasing size. 
	The table reports, for each working precision, 
	the percentage of matrix multiplications done in precision $u^{1/2}$ or lower ($M_{low}$), the total execution
	time in seconds for Algorithm~\ref{alg.mpps} ($T_{tot}$)
	and its fixed-precision counterpart ($T_{fix}$), and 
	the percentage of time spent calculating the required matrix powers in precision $u$ ($T_{pow}$), estimating the matrix norm ($T_{est}$), executing the Horner's scheme with matrix multiplications and additions in mixed precisions ($T_{hon}$), 
	and assembling the coefficient matrices ($T_{coe}$).}
	\resizebox{1.01\linewidth}{!}{
		\input{matlab/tabs/table_profile_0064_0256.tex}
	}\label{tab:profiling}
\end{table}

Table~\ref{tab:profiling} compares the execution time of our MATLAB
implementations of Algorithm~\ref{alg.mpps} and its fixed-precision counterpart from~\cite{fahi19}, run in 64 and 256 decimal significant digits, on the three classes of matrices (which, as mentioned before, are prior to the scaling)
\begin{verbatim}
	A = gallery('smoke', n);		   % contians 2*n nonzero elements
	B = 100 * triu(randn(n),1);  % strictly upper-triangular
	C = gallery('lotkin', n);		  % full
\end{verbatim}
where $n$ varies between 20 and 1000. A matrix generated by \texttt{A} is a perturbed bidiagonal matrix with $1$ in its $(n,1)$ position, making it increasingly sparse as the dimension $n$ increases.
The second class \texttt{B} consists of strictly upper triangular matrices, so these matrices are nilpotent and the positive integer powers of the matrices are in general increasingly sparse. The third class \texttt{C} contains full matrices with eigenvalues exponentially decaying in modulus.

The columns of $M_{low}$ in Table~\ref{tab:profiling} show that typically at least 20\% of the matrix multiplications were performed in precision $u^{1/2}$ or lower, and this is consistent with our observation in Table~\ref{tab.exp_cauchy}. 
For all three classes of matrices, the norm estimation ($T_{est}$) of the associated coefficient matrices $B_i(\cdot)$, $i=0\colon r$ and the $s$th power is relatively expensive for small matrices, but its overhead has a lesser impact on the overall runtime as the size of the matrices increases; and it is typically negligible when the matrix size reaches 200.

For the matrices in \texttt{A}, Algorithm~\ref{alg.mpps} is slower than its fixed-precision counterpart in all tested cases, this is due to the high sparsity of these matrices that the operations of matrix scaling and addition performed in assembling the coefficient matrices $B_i(\cdot)$ ($T_{coe}$) become dominant in cost.
For the matrices in \texttt{B} and \texttt{C}, the calculation of the matrix powers in the working precision $u$ ($T_{pow}$) becomes the most expensive operation for matrices of size larger than 200 and 100, respectively;
on the other hand, the operation costs measured by $T_{coe}$ only depends quadratically on the $n$, so this time ratio generally 
declines as $n$ increases.
We start to see the advantage of Algorithm~\ref{alg.mpps} in terms of 
total execution time for matrices in \texttt{C} of size at least 500 when the working precision is 64 decimal digits, and the size threshold is lowered to 100 when the working precision increases to 256 decimal digits; in the working precision of 64 decimal digits the advantage is marginal, but in the higher working precision Algorithm~\ref{alg.mpps} saves 17.4\% of runtime for matrices of size 1000.
Indeed, the actual performance gain is less than what is predicted by the value of $C_r$ in~\eqref{eq:Cr} via counting the theoretical numerical operations, but we do see a trend that the advantage of the algorithm is enlarged as the matrix size increases or the working precision becomes higher.

\begin{figure}
		\begin{subfigure}{.5\linewidth}
			\centering
			\includegraphics[height=4.3cm]{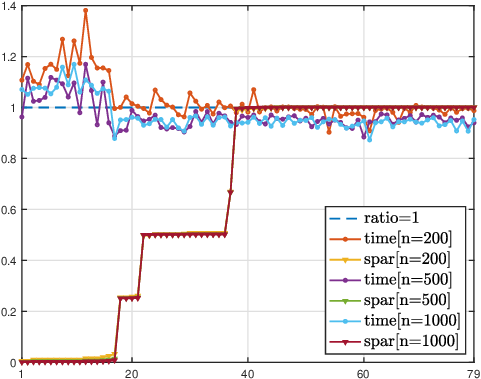}
			\caption{$u = 10^{-64}$.}
		\end{subfigure}
		\hspace{-25pt}
		\begin{subfigure}{.5\linewidth}
			\centering
			\includegraphics[height=4.3cm]{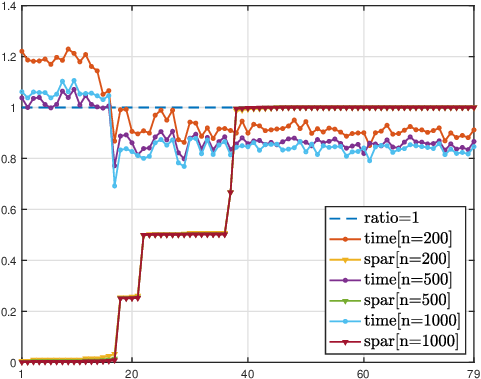}
			\caption{$u = 10^{-256}$.}
		\end{subfigure}
				\vspace{-10pt}
				
		\caption{The ratio between the runtime of Algorithm~\ref{alg.mpps} and the runtime of its fixed-precision counterpart from~\cite{fahi19} on matrices of different sizes, denoted by $\texttt{time[size]}$. The matrices are sorted in ascending order of their sparsity measured by the proportion of nonzero elements, characterized by the curves $\texttt{spar[size]}$.}
		\label{fig.exp_taylor_time}
\end{figure}

Finally, from the test set we used in Section~\ref{sect.numerical.accuracy.speedup} we take the 79 test matrices that have variable dimension parametrized by $n$, and we compare the runtime of Algorithm~\ref{alg.mpps} and its fixed-precision counterpart
on these matrices with different choices of $n$; the results are presented in Figure~\ref{fig.exp_taylor_time}.
Not surprisingly, we see Algorithm~\ref{alg.mpps} cannot gain any advantage in runtime on matrices that are very sparse (with at most 20\% of nonzero matrix elements); in the working precision of 64 decimal digits the mixed-precision algorithm is 5\% to 10\% faster on the rest of matrices with size larger than 200. 
In the higher precision of 256
decimal digits, we can already see recognizable acceleration for matrices of size 200, and this speed-up is typically around 20\% as the matrix size increases to 1000.

\begin{algorithm2e}[t]
	\setstretch{1}
	\caption{Mixed-precision Paterson--Stockmeyer scheme for matrix polynomials with scalar coefficients.}
	\label{alg.mpps.gen}
	\SetAlgoVlined
	{\nonl
		\begin{minipage}{0.95\linewidth}
			\pretolerance=100\tolerance=200\hyphenpenalty=5
			Given $X\in\mathbb{C}^{\nbyn}$ and a set of polynomial coefficients $\{b_i\}_{i=0}^m$, this algorithm computes the matrix polynomial with scalar coefficients decaying in modulus
			$P\equiv p_m(X)$ in the form of~\eqref{eq:ps} using the Paterson--Stockmeyer scheme in floating-point arithmetic. The algorithm starts with the user-specified precision $u$ and potentially uses multiple lower precisions $u_i\ge u$
			aiming to produce a relative error of order $nu$.
	\end{minipage}}
	\BlankLine
	$s \gets \ceil{\sqrt{m}}$ \; 
	$r \gets \floor{m/s}$ \;
	$\nu \gets r+1$ \;
	$u_0 \gets u$ \;
	$B_0 \gets b_0I$ in precision $u_0$\; 
	\For{$j \gets 1$ \algto $s-1$}{
		$\X_{j} \gets \X_{j-1} X$ in precision $u_0$ \;
		$B_0 \gets B_0 + b_j\X_{j}$ in precision $u_0$ \; 
	}
	$Y \gets \X_{s-1} X$ in precision $u_0$ \;
	\For{$i \gets 1$ \algto $r$}{
		Assemble $B_i$ using elements in $\X\cup\{I\}$ and then estimate $\normi{B_i}$\;
		$u_i \gets \normi{B_{i-1}}u_{i-1}/(\normi{B_i}\normi{Y}^i)$\;
	}
	$\nu\gets \min\{i\colon u_i\ge \delta u\}$, $u_{\nu-1},u_{\nu-2},\dots,u_{1}\gets u$ \;
	$P = B_r$ \;
	\For{$i \gets r$ \algdownto $1$}{
		Compute $P\gets PY$ in precision $u_i$\;
		Form $P\gets P + B_{i-1}$ in precision $u_{i-1}$ \;	}
	\Return{$P$}\;
	\SetAlgoNoLine
\end{algorithm2e}

\section{A mixed-precision Paterson--Stockmeyer algorithm for general polynomials of matrices}\label{sect.general.framework}
Following the discussion in Section~\ref{sect.mpps.exp.taylor}, we now exploit the mixed-precision 
PS framework for the computation of general matrix polynomials with scalar coefficients decaying in modulus.
The algorithm is presented as~Algorithm~\ref{alg.mpps.gen}. 

We will test Algorithm~\ref{alg.mpps.gen} in the next subsections on different types of polynomials of matrices arising from the computation of exponential-like matrix functions.

\begin{figure}
	\begin{subfigure}{1\linewidth}
		\centering
		\includegraphics[height=4.3cm]{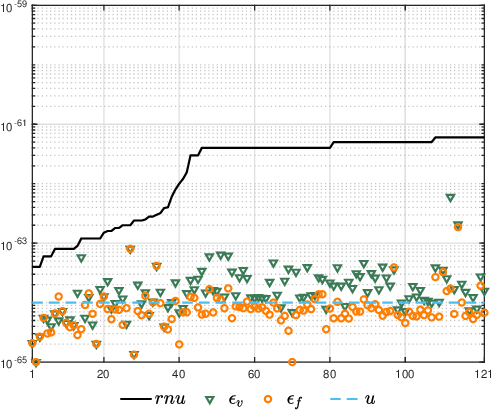} 
		\includegraphics[height=4.3cm]{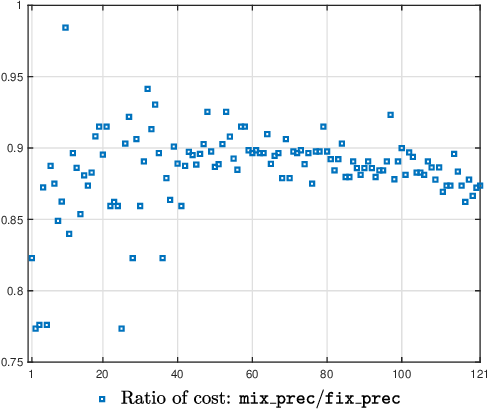}
		\caption{$p_m=p_{mm}$ from the $[m/m]$ Pad\'e approximant of the matrix exponential.}
	\end{subfigure} \\
	\vspace{-15pt}
	
	\begin{subfigure}{1\linewidth}
		\centering
		\includegraphics[height=4.3cm]{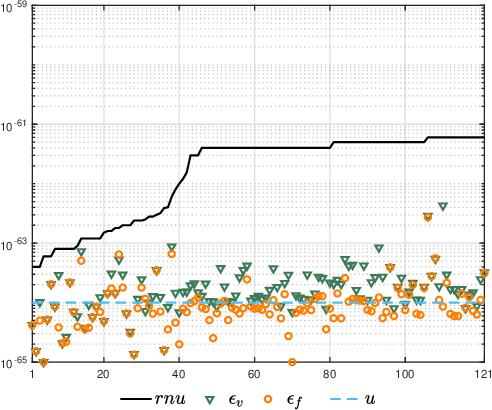} 
		\includegraphics[height=4.3cm]{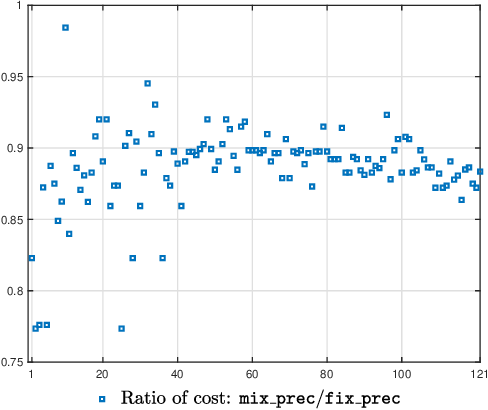}
		\caption{$p_m=q_{mm}$ from the $[m/m]$ Pad\'e approximant of the matrix exponential.}
	\end{subfigure} \\
	\vspace{-15pt}
	
	\begin{subfigure}{1\linewidth}
		\centering
		\includegraphics[height=4.3cm]{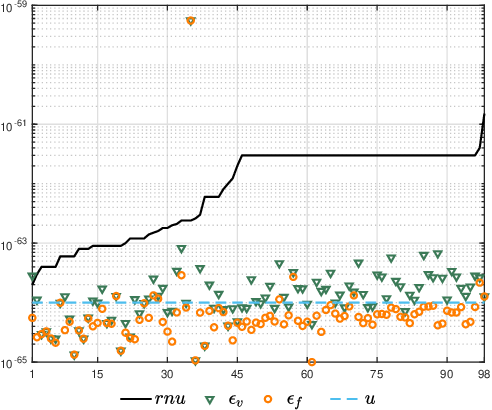}
		\includegraphics[height=4.3cm]{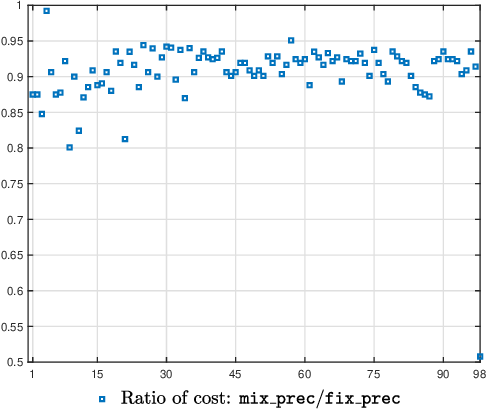}
		\caption{$p_m=c_{m}$ from the Taylor approximant of the matrix cosine.}
	\end{subfigure}
	\vspace{-15pt}
	
	\caption{Left: The relative error $\epsilon_v=\normi{\widehat{p}_m-{p}_m(X)}/\normi{{p}_m(X)}$ produced by Algorithm~\ref{alg.mpps.gen} in precision $u = 10^{-64}$ compared with the relative error $\epsilon_f$ produced by the fixed-precision counterpart from~\cite{ahl22} or~\cite{fahi19} with $s=\ceil{\sqrt{m}}$ on various matrices with $2\le n\le 100$. Right: The associated approximate ratio of cost $C_r$ in~\eqref{eq:Cr}.}
	\label{fig.exp_pade_cos_taylor_ap}
\end{figure}

\subsection{Pad\'e approximants to the matrix exponential}\label{sect.pade.exp}
The numerator and denominator of the Pad\'e approximant for the matrix exponential are computed separately 
by means of the PS method in~\cite{fahi19}. Therefore, similar to the experiments carried out in Section~\ref{sect.numerical.exp.taylor}, we now feed the working precision $u$ and the matrix from the same test set used in Section~\ref{sect.numerical.accuracy.speedup} 
to the Pad\'e-based scaling and squaring algorithm~\cite[Alg.~4.1]{fahi19}, in order to extract the test data, namely, 
the matrix polynomials $p_{km}(X)$ and $q_{km}(X)$ from the $[k/m]$ Pad\'e approximants $r_{km}(X)= q_{km}^{-1}(X)p_{km}(X)$, where $p_{km}(X)$ and $q_{km}(X)$ are of degrees at most $k$ and $m$, respectively.
We then compare Algorithm~\ref{alg.mpps.gen} with Algorithm~\ref{alg.fpps} 
(both specified to the numerator and denominator Pad\'e approximants to the matrix exponential, respectively)
in precision $u$. 

Figure~\ref{fig.exp_pade_cos_taylor_ap}~$(a)$--$(b)$ reports the result. We observe that the relative errors produced by Algorithm~\ref{alg.mpps.gen} have the same order of magnitude as those produced by the fixed-precision counterpart. The reduction in computational complexity of the algorithm is on average around 10\%, which is only half the reduction obtained when applying the algorithm to the Taylor approximants of the matrix exponential (cf.~Figure~\ref{fig.exp_taylor_ap}~$(a)$). 

The smaller (approximately halved) chosen degrees $m$ in the matrix Pad\'e approximants than the matrix Taylor approximant have prevented the algorithm from carrying out more low-precision matrix multiplications, outweighing the fact that the algorithm can potentially employ lower precisions early (since the scalar coefficients of the matrix polynomials from the Pad\'e approximant decay much faster in magnitude).

\subsubsection{Low-precision variant}
We also implemented a low-precision-oriented variant of
Algorithm~\ref{alg.mpps.gen} that selects precisions from the IEEE double, IEEE single, and bfloat16 half precisions, aiming to achieve accuracy of level of the unit roundoff of double precision. 
This combination of precisions has attracted a lot of interest in scientific computing, especially after the half precision is supported on actual hardware~\cite{hth19}.
In the variant algorithm, the $u_i$ calculated by~\eqref{eq:choose-ui} will be set to the unit roundoff of the nearest higher precision among the three available precisions, that is, 
$u_i \gets \max\{u^* \le u_i \colon u^*\in \{\uf,\us, u\} \}$,
where
$\uf
\approx 3.9\times10^{-3}$, $\us
\approx 6.0\times10^{-8}$, and $u
\approx 1.1\times 10^{-16}$.
Here bfloat16 instead of the IEEE half precision is used because the latter has the drawback of having a limited range, the smallest representable positive number being $x_{\min}\approx5.96\times 10^{-8}$, while the former has 
much larger exponential range with $x_{\min}\approx1.18\times 10^{-38}$ and this makes loss of accuracy due to underflow less likely to happen; this is particularly important for Algorithm~\ref{alg.mpps.gen} and the variant, whose spirit is computing smaller numbers in a lower precision.

We tested the low-precision variant of Algorithm~\ref{alg.mpps.gen} on the same set of $121$ matrices used in
Section~\ref{sect.numerical.accuracy.speedup}, in which case the degree $m$ of the \pade polynomial and the scaling parameter $\ell$ were chosen by the double-precision-oriented scaling and squaring algorithm~\cite{alhi09a} (cf.~Sect.~\ref{sect.numerical.exp.taylor}).
We found the results are similar to that reported in Figure~\ref{fig.exp_pade_cos_taylor_ap}~$(a)$--$(b)$ in terms of accuracy, while the theoretical efficiency gain of the mixed-precision algorithm was further reduced, at a typical level of less than 10\%: in many cases the algorithm just degenerated to solely choose double precision.

The behaviour of the variant algorithm is expected since Algorithm~\ref{alg.mpps.gen} is mainly designed for and useful in high precision computing, especially in an arbitrary precision environment, as we have seen in the experiments of Section~\ref{sect.numerical.exp.taylor}.
This is because the chosen degree $m$ is generally small in a low working precision, and this restricts more matrix products being formed in lower precisions and hence limits the efficiency gain of the mixed-precision algorithm. Also, in the double-single-half precision environment the algorithm loses the complete freedom to choose arbitrary precision and it in many cases has to use an unnecessarily higher precision for accuracy.

\subsection{Taylor approximants to the matrix cosine}
Finally, we test Algorithm~\ref{alg.mpps.gen} using the Taylor approximants of the matrix cosine, denoted by $c_m$. We use the $98$ nonnormal test matrices from~\cite{ahl22} (with the $2\times 2$ nilpotent matrix therein excluded) that have size between $4$ and $100$. The degree $m$ and the scaling parameter $\ell$ for these matrices are determined by the 
Taylor-based scaling and recovering algorithm of~\cite{ahl22}.

Figure~\ref{fig.exp_pade_cos_taylor_ap}~$(c)$ shows that the accuracy of Algorithm~\ref{alg.mpps.gen} in this case is similar to that of the fixed-precision counterpart (Algorithm~\ref{alg.fpps} specified to the Taylor approximant of the matrix cosine).
The moduli of the coefficients in the Taylor approximants of the matrix cosine also decay considerably faster than that of the matrix exponential, and for the same reasons as discussed in 
Section~\ref{sect.pade.exp},
the theoretical reduction in computational cost of the algorithm is marginal, at typically less than 10\%. 
Moreover, in the algorithm we are actually evaluating the polynomial $c_m(X^2)$ via the PS scheme since the polynomial only contains even powers of $X$, and this makes the chosen degree $m$ even smaller.

\section{Conclusions}\label{sect.concl}
In this work we have developed a mixed-precision Paterson--Stockmeyer (PS) method for evaluating matrix polynomials with scalar coefficients that decay in modulus. The key idea is to perform computations on data of small magnitude in low precision, and driven by this idea we show that, if the coefficients satisfy a certain decaying property, then we can exploit a set of suitably chosen lower precisions (relative to the working precision) in the evaluation of the polynomial via a Horner scheme and still achieve the same level of accuracy as the evaluation being done solely in the working precision.

We applied the method to the computation of the Taylor approximants of the matrix exponential and showed  
the applicability of the mixed-precision PS framework to the existing scaling and squaring algorithms for the matrix exponential.
The algorithm (Algorithm~\ref{alg.mpps}) switches to lower precisions when appropriate, and it is readily employable by the Taylor-based scaling and squaring algorithms for the matrix exponential.
Taking advantage of the generality of the mixed-precision PS framework, we also propose an algorithm (Algorithm~\ref{alg.mpps.gen}) for the computation of general polynomials of matrices and demonstrated its efficiency on polynomials from the Pad\'e approximant of the matrix exponential and the Taylor approximant of the matrix cosine.

Numerical experiments show comparable accuracy of our mixed-precision PS 
algorithms to its fixed-precision counterparts (in the working precision). 
By gauging the computational complexity in terms of the number of matrix products and assuming the complexity is linearly proportional to the number of used digits, we obtain the theoretical computational savings of the algorithms in various settings.
We also compared the actual runtime of the algorithms in software and 
it confirms the performance gain of the new algorithm; while it is  
less than what is predicted by counting the theoretical numerical operations, the trend is clear: the advantage of the algorithm is enlarged as the matrix size increases or the working precision becomes higher. 
For matrices of larger size the extra $O(n^2)$ flops (compared with the fixed-precision PS method) becomes more negligible. 
On the other hand, the number of required matrix
multiplications in the working precision can be reduced from $m$ to approximately $\sqrt{m}$ and this reduction is substantial for large degree $m$, which is often the case for arbitrary precision environment, where the working precision can be much higher than the IEEE double precision.

As a possible direction of future research, the proposed mixed-precision PS method may be improved by using mixed precisions also in 
computing the first $s$ positive matrix powers: the higher powers can be computed in a lower precision and 
the accuracy loss will be mitigated when these powers are multiplied by the scalar coefficients that decay in modulus; meanwhile, this mixed-precision scheme might not compute the $s$th power $X^s$ to the required accuracy, and in this case binary powering in the working precision may be used for the purpose.
This idea is probably more relevant at the setting where the working precision is high, in which case the chosen degree $m$ is high so that $s$ is large, giving more potential benefits for the use of mixed precisions.

\appendix
\section{Proof of Theorem~\ref{thm:ps-mp-err-norm}}\label{append:proof-ps-mp-err}
The proof is by induction on the normwise error of the computed $\varphi_{j}$, which is defined as
\begin{equation}\label{eq:induction-object}
	E_{r-i}:=\wh{\varphi}_{r-i}-\varphi_{r-i},\quad i=0\colon r-\nu+1.
\end{equation}
From the hypotheses of theorem, we have
$\norm{E_r}=\norm{\wh{\varphi}_r-\varphi_r} = \norm{\wh{B}_r-B_r} \le u_r \norm{B_r}$.
Consider 
\begin{align*}
	E_{r-1} &= \wh{\varphi}_{r-1}-\varphi_{r-1} = 
	\fl_{r-1}\big(\wh{B}_{r-1}+\fl_{r}(\wh{Y}\wh{\varphi}_{r})\big) - B_{r-1} - Y\varphi_{r} \\
	&:=  E_{r-1,s} + E_{r-1,p} + E_{r-1,a},
\end{align*}
where $E_{r-1,s}$ denotes the error in the computed matrix coefficient  $\wh{B}_{r-1}$,
$E_{r-1,p}$ represents the error occurring in forming the matrix product, and 
$E_{r-1,a}$ accounts for the error in the subsequent matrix addition. We have
\begin{equation*}
	\norm{E_{r-1,s}} \le \gamma_1^{r-1}\norm{B_{r-1}},\quad 
	\norm{E_{r-1,p}} \le \gamma_n^r\norm{Y}\norm{B_r}, 
\end{equation*}
and 
\begin{equation*}
	\norm{E_{r-1,a}}\le \gamma_1^{r-1}\norm{B_{r-1} + E_{r-1,s} + YB_r +E_{r-1,p}}.
\end{equation*}
So by using~\eqref{eq:gamma-theta} we have
\begin{align*}
	\norm{E_{r-1}} &\le
	\left(\gamma_n^r + \gamma_1^{r-1}+\gamma_n^r\gamma_1^{r-1}\right)\norm{Y}\norm{B_r} + \gamma_2^{r-1}\norm{B_{r-1}} \\
	&\le \gamma_{n\theta_{r,r-1} +1}^{r-1} \norm{Y}\norm{B_r} + \gamma_2^{r-1}\norm{B_{r-1}}.
\end{align*}
Then we have, for the \textit{base case} with $i=2$ in~\eqref{eq:induction-object},
\begin{align*}
	E_{r-2} &= \wh{\varphi}_{r-2}-\varphi_{r-2} = 
	\fl_{r-2}\big(\wh{B}_{r-2}+\fl_{r-1}(\wh{Y}(\varphi_{r-1}+E_{r-1}))\big) - B_{r-2} - Y\varphi_{r-1} \\
	&=  YE_{r-1} + E_{r-2,s} + E_{r-2,p} + E_{r-2,a},
\end{align*}
where, with similar notations,
\begin{equation*}
	\norm{E_{r-2,s}} \le \gamma_1^{r-2}\norm{B_{r-2}},\quad 
	\norm{E_{r-2,p}} \le \gamma_n^{r-1}\norm{Y}\norm{\varphi_{r-1}+E_{r-1}}, 
\end{equation*}
and 
\begin{equation*}
	\norm{E_{r-2,a}}\le \gamma_1^{r-2}\norm{B_{r-2} + E_{r-2,s} + Y\left(\varphi_{r-1}+E_{r-1}\right) +E_{r-2,p}},
\end{equation*}
so we have
\begin{align*}
	\norm{E_{r-2}} &\le \gamma_2^{r-2}\norm{B_{r-2}} + \big(1+\gamma_{n\theta_{r-1}+1}^{r-2}\big)\norm{Y}\norm{E_{r-1}} +  
	\gamma_{n\theta_{r-1}+1}^{r-2}\norm{Y}\norm{\varphi_{r-1}} \\
	& \le
	\gamma_{f_{r-2,r}}^{r-2}\norm{Y}^2\norm{B_r} + \gamma_{f_{r-2,r-1}}^{r-2}\norm{Y}\norm{B_{r-1}} + \gamma_{f_{r-2,r-2}}^{r-2}\norm{B_{r-2}},
\end{align*}
where 
\begin{align*}
	f_{r-2,r} &= n\theta_{r,r-2}+(n+1)\theta_{r-1,r-2}+1 \\
	f_{r-2,r-1} & = (n+2)\theta_{r-1,r-2}+1 \\
	f_{r-2,r-2} & = 2.
\end{align*}
For the \textit{inductive step}, now assume a bound on~\eqref{eq:induction-object} with $i=k$ of the following form:
\begin{equation}\label{eq:ps-mp-e_r-k}
	\norm{E_{r-k}}  \le
	\gamma_{f_{r-k,r}}^{r-k}\norm{Y}^k\norm{B_r} + \gamma_{f_{r-k,r-1}}^{r-k}\norm{Y}^{k-1}\norm{B_{r-1}} + \cdots +  \gamma_{f_{r-k,r-k}}^{r-k}\norm{B_{r-k}},
\end{equation}
where
\begin{align*} 
	 f_{r-k,r} &= n\theta_{r,r-k} + (n+1)\left(\theta_{r-1,r-k}+\theta_{r-2,r-k}+\dots+\theta_{r-k+1,r-k}\right) + 1,   \\
	 f_{r-k,r-1} &= (n+2)\theta_{r-1,r-k} + (n+1)\left(\theta_{r-2,r-k}+\theta_{r-3,r-k}+\dots+\theta_{r-k+1,r-k}\right) + 1,   \\
				&\vdots  \\
	f_{r-k,r-k+1} &= (n+2)\theta_{r-k+1,r-k}+1, \\
	f_{r-k,r-k} &= 2.
\end{align*}
Then we have
\begin{align*}
	E_{r-(k+1)} &= \wh{\varphi}_{r-(k+1)}-\varphi_{r-(k+1)} \\
	&= 
	\fl_{r-(k+1)}\big(\wh{B}_{r-(k+1)}+\fl_{r-k}(\wh{Y}(\varphi_{r-k}+E_{r-k}))\big) - B_{r-(k+1)} - Y\varphi_{r-k} \\
	&=  YE_{r-k} + E_{r-(k+1),s} + E_{r-(k+1),p} + E_{r-(k+1),a},
\end{align*}
where
\begin{equation*}
	\left\|E_{r-(k+1),s}\right\| \le \gamma_1^{r-(k+1)}\left\|B_{r-(k+1)}\right\|,\quad 
	\left\|E_{r-(k+1),p}\right\| \le \gamma_n^{r-k}\norm{Y}\norm{\varphi_{r-k}+E_{r-k}}, 
\end{equation*}
and 
\begin{equation*}
	\left\|E_{r-(k+1),a}\right\|\le \gamma_1^{r-(k+1)}\left\|B_{r-(k+1)} + E_{r-(k+1),s} + Y\left(\varphi_{r-k}+E_{r-k}\right) +E_{r-(k+1),p}\right\|,
\end{equation*}
and we have
\begin{multline*}
	\left\|E_{r-(k+1)}\right\|  \le
	\gamma_2^{r-(k+1)}\left\|B_{r-(k+1)}\right\| + \big(1+\gamma_{n\theta_{r-k,r-(k+1)}+1}^{r-(k+1)}\big)\norm{Y}\norm{E_{r-k}} 
	\\ +\gamma_{n\theta_{r-k,r-(k+1)}+1}^{r-(k+1)}\norm{Y}\norm{\varphi_{r-k}},
\end{multline*}
where $\varphi_{r-k}=Y^kB_r + Y^{k-1}B_{r-1} + \cdots 
+ B_{r-k}$. Then by writing the gamma constants in~\eqref{eq:ps-mp-e_r-k} as
\begin{equation*}
	\gamma_{f_{r-k,j}}^{r-k} = \gamma_{\theta_{r-k,r-(k+1)}f_{r-k,j}}^{r-k+1},\quad j=r-k \colon r,
\end{equation*}
and by noting the bounds  
\begin{equation*}
	\gamma_{\theta_{r-k,r-(k+1)}f_{r-k,j}}^{r-k+1}\big(1+\gamma_{n\theta_{r-k,r-(k+1)}+1}^{r-(k+1)}\big)
	+ \gamma_{n\theta_{r-k,r-(k+1)}+1}^{r-(k+1)} \le  \gamma_{\theta_{r-k,r-(k+1)}(n+f_{r-k,j})+1}^{r-k+1},
\end{equation*}
we have
\begin{multline*}
	\left\|E_{r-(k+1)}\right\|  \le
	\gamma_{f_{r-(k+1),r}}^{r-(k+1)}\norm{Y}^{k+1}\norm{B_r} + \gamma_{f_{r-(k+1),r-1}}^{r-(k+1)}\norm{Y}^{k}\norm{B_{r-1}} + \cdots 
	\\ + \gamma_{f_{r-(k+1),r-(k+1)}}^{r-(k+1)}\left\|B_{r-(k+1)}\right\|,
\end{multline*}
where $f_{r-(k+1),r-(k+1)} = 2$ and
\begin{align*} 
	f_{r-(k+1),r} &= \theta_{r-k,r-(k+1)}\left(n+f_{r-k,r}\right)+1 \\
	 &= \begin{aligned}[t]	
	 	&n\theta_{r,r-(k+1)} + \\
	 	&(n+1)\left(\theta_{r-1,r-(k+1)}+
	 \theta_{r-2,r-(k+1)}+
\dots+\theta_{r-k,r-(k+1)}\right) + 1, \end{aligned}
\end{align*}	
and similarly, $f_{r-(k+1),j}= \theta_{r-k,r-(k+1)}\left(n+f_{r-k,j}\right)+1$ for $j=r-k\colon r-1$, so
\begin{align*} 	
		 f_{r-(k+1),r-1} &= \begin{aligned}[t]	
		 	&(n+2)\theta_{r-1,r-(k+1)} + \\
		 	&(n+1)\left(\theta_{r-2,r-(k+1)} +\theta_{r-3,r-(k+1)}+ \cdots  
		 	+\theta_{r-k,r-(k+1)}\right) + 1,   
		 \end{aligned}\\
	&\vdots  \\
	f_{r-(k+1),r-k+1} &= (n+2)\theta_{r-k+1,r-(k+1)}+ (n+1)\theta_{r-k,r-(k+1)} +1, \\
	f_{r-(k+1),r-k} &= (n+2)\theta_{r-k,r-(k+1)} + 1.
\end{align*}
This proves~\eqref{eq:induction-object} for $i=k+1$ and so the proof is completed by induction. 

\section*{Acknowledgments}
This work was carried out when the author was a research associate in the Numerical Linear Algebra Group at The University of Manchester.
The author is grateful to Prof.~Nicholas J. Higham for financial support through his grant from the Royal Society and for valuable comments on early versions of the manuscript, and to Dr.~Massimiliano Fasi for his useful suggestions.
The author also thanks the anonymous referees for their comments, which helped him improve the experimental sections and the presentation of the paper.

\bibliographystyle{siamplain}
\bibliography{strings,njhigham,notebib}

\end{document}

%% file: matlab/tabs/table_cauchy_100.tex
\begin{tabularx}{\textwidth}{@{\extracolsep{\fill}}crrrrr}
\toprule
decimal digits & $m$ & $s$ & $r$ & $(d_1,d_2,\dots,d_r)$  & $C_r$ \\
\midrule 
32 & 42 & 7 & 6 & $(30, 25, 18, 11, 3, 1)$ & 27.1\% \\
64 & 64 & 8 & 8 & $(61, 55, 47, 38, 28, 18, 7, 1)$ & 26.8\% \\
128 & 100 & 10 & 10 & $(124, 115, 104, 92, 78, 64, 49, 34, 18, 1)$ & 24.7\% \\
256 & 182 & 14 & 13 & $(248, 234, 217, 197, 176, 154, 131, 107, 82, 57, 31, 4, 1)$ & 25.4\% \\
\bottomrule
\end{tabularx}

%% file: matlab/tabs/table_profile_0064_0256.tex
\begin{tabularx}{1.095\textwidth}{@{\extracolsep{\fill}}rr|rrrrrrr|rrrrrrr}
\toprule
\multicolumn{2}{c|}{} & \multicolumn{7}{c|}{ 64 decimal digits} & \multicolumn{7}{c}{ 256 decimal digits} \\
 & $n$ & $M_{low}$ & $T_{pow}$ & $T_{est}$ & $T_{hon}$ & $T_{coe}$ & $T_{tot}$ & $T_{fix}$ & $M_{low}$ & $T_{pow}$ & $T_{est}$ & $T_{hon}$ & $T_{coe}$ & $T_{tot}$ & $T_{fix}$\\
\midrule
\verb|A| &  20 & 27\% &  4\% & 52\% & 16\% & 29\% &  0.0 &  0.0 & 24\% &  5\% & 33\% & 18\% & 44\% &  0.1 &  0.0 \\
         &  50 & 27\% &  7\% & 29\% & 20\% & 43\% &  0.1 &  0.0 & 24\% &  8\% & 16\% & 23\% & 54\% &  0.2 &  0.1 \\
         & 100 & 27\% &  9\% & 21\% & 19\% & 51\% &  0.2 &  0.1 & 24\% & 19\% &  7\% & 15\% & 59\% &  0.9 &  0.7 \\
         & 200 & 27\% & 15\% & 10\% & 14\% & 60\% &  0.9 &  0.7 & 24\% & 17\% &  4\% & 13\% & 66\% &  4.0 &  3.3 \\
         & 500 & 27\% & 11\% &  6\% & 12\% & 71\% & 10.1 &  9.4 & 24\% & 10\% &  3\% & 10\% & 77\% & 46.8 & 44.2 \\
         & 1000 & 27\% & 10\% &  6\% & 16\% & 68\% & 44.7 & 38.2 & 24\% & 10\% &  3\% & 10\% & 77\% & 191.4 & 172.7 \\
\midrule
\verb|B| &  20 & 24\% &  4\% & 62\% &  9\% & 25\% &  0.0 &  0.0 & 23\% &  5\% & 52\% &  8\% & 35\% &  0.0 &  0.0 \\
         &  50 & 22\% & 14\% & 35\% & 13\% & 38\% &  0.1 &  0.0 & 25\% & 18\% & 21\% & 10\% & 51\% &  0.1 &  0.1 \\
         & 100 & 24\% & 24\% & 16\% & 21\% & 38\% &  0.2 &  0.1 & 23\% & 30\% &  8\% & 15\% & 48\% &  0.7 &  0.6 \\
         & 200 & 22\% & 38\% &  7\% & 26\% & 30\% &  0.9 &  0.8 & 21\% & 40\% &  3\% & 17\% & 40\% &  4.3 &  4.1 \\
         & 500 & 25\% & 43\% &  3\% & 28\% & 26\% & 12.7 & 12.5 & 20\% & 38\% &  1\% & 19\% & 42\% & 66.4 & 72.6 \\
         & 1000 & 22\% & 42\% &  2\% & 31\% & 25\% & 87.2 & 89.1 & 24\% & 46\% &  1\% & 26\% & 27\% & 381.5 & 433.2 \\
\midrule
\verb|C| &  20 & 25\% &  7\% & 52\% & 13\% & 28\% &  0.0 &  0.0 & 26\% & 10\% & 37\% & 12\% & 40\% &  0.0 &  0.0 \\
         &  50 & 25\% & 24\% & 23\% & 22\% & 31\% &  0.1 &  0.0 & 22\% & 28\% & 12\% & 24\% & 36\% &  0.2 &  0.2 \\
         & 100 & 27\% & 34\% & 11\% & 31\% & 24\% &  0.2 &  0.2 & 24\% & 42\% &  4\% & 26\% & 28\% &  0.8 &  0.9 \\
         & 200 & 27\% & 42\% &  4\% & 37\% & 16\% &  1.3 &  1.3 & 24\% & 48\% &  2\% & 28\% & 23\% &  6.1 &  6.8 \\
         & 500 & 27\% & 47\% &  1\% & 40\% & 12\% & 17.6 & 18.4 & 23\% & 49\% &  1\% & 27\% & 23\% & 93.1 & 109.8 \\
         & 1000 & 27\% & 46\% &  1\% & 41\% & 12\% & 137.7 & 145.5 & 23\% & 55\% &  0\% & 31\% & 13\% & 636.6 & 770.8 \\
\bottomrule
\end{tabularx}